\pdfoutput=1
\documentclass[11pt,oneside,reqno]{amsart}
\usepackage{amsmath, amssymb, amsthm, amsfonts, mathrsfs}

\setbox0=\hbox{$+$}
\newdimen\plusheight
\plusheight=\ht0
\def\+{\;\lower\plusheight\hbox{$+$}\;}

\setbox0=\hbox{$-$}
\newdimen\minusheight
\minusheight=\ht0
\def\-{\;\lower\minusheight\hbox{$-$}\;}

\setbox0=\hbox{$\cdots$}
\newdimen\cdotsheight
\cdotsheight=\plusheight
\def\cds{\lower\cdotsheight\hbox{$\cdots$}}

\newtheorem{question}{Question}
\newtheorem{conjecture}{Conjecture}

\theoremstyle{definition}

\usepackage{graphicx, epsfig}
\theoremstyle{definition}


\numberwithin{equation}{section}
\theoremstyle{plain}
\newtheorem{theorem}{Theorem}[section]

\newtheorem{lemma}[theorem]{Lemma}
\newtheorem{claim}[theorem]{Claim}
\usepackage{titlesec}

\usepackage{graphicx, epsfig}
\theoremstyle{definition}
\usepackage{graphicx} 
\usepackage{color} 
\usepackage{transparent} 
\titleformat*{\section}{\LARGE\bfseries}

\makeatletter
\renewcommand\section{\@startsection{section}{1}{\z@}%
                                  {-3.5ex \@plus -1ex \@minus -.2ex}%
                                  {2.3ex \@plus.2ex}%
                                  {\normalfont\large\bfseries}}
\makeatother
\usepackage{graphicx} 
\usepackage{color} 
\usepackage{transparent} 
\usepackage{geometry}
\begin{document}
\title{Realizing Algebraic Invariants of Hyperbolic Surfaces}

\author{BoGwang Jeon}
\begin{abstract}
Let $S_g$ ($g\geq 2$) be a closed surface of genus $g$. Let $K$ be any real number field and $A$ be any quaternion algebra over $K$ such that $A\otimes_K\mathbb{R}\cong M_2(\mathbb{R})$. 
We show that there exists a hyperbolic structure on $S_g$ such that $K$ and $A$ arise as its invariant trace field and invariant quaternion algebra.
\end{abstract} 
\maketitle

\section{Introduction}
The invariant trace field and the quaternion algebra are basic algebraic invariants of a Kleinian group. 
Questions around these invariants have inspired extensive research in the study of hyperbolic $3$-manifolds. 
The following Realization Conjecture due to W. Neumann \cite{walter} is one of the fundamental questions along these lines:

\begin{conjecture}
Let $K$ be any non-real complex number field and $A$ be any quaternion algebra over $K$. Then there exists a hyperbolic $3$-manifold $M$ such that $K$ and $A$ arise as the invariant trace field and quaternion algebra of $M$.
\end{conjecture}

For the $2$-dimensional case, since Mostow's rigidity theorem does not hold, trace fields and quaternion algebras are no longer topological invariants. But, in this case, we can ask the following analogous question instead, as suggested in \cite{walter}:
\begin{question}\label{qe1}
Let $S_g$ ($g \geq 2$) be a closed surface of genus $g$, $K$ be any real number field and $A$ be any quaternion algebra over $K$ such that $A\otimes_K\mathbb{R}\cong M_2(\mathbb{R})$.\footnote{We say $K$ is a real number field if $K\subset \mathbb{R}$ and $[K:\mathbb{Q}]<\infty$. 
Note that the condition on $A$ (i.e. $A\otimes_K\mathbb{R}\cong M_2(\mathbb{R})$) is necessary (see Section \ref{qa}).} 
Is there a hyperbolic structure on $S_g$ whose invariant trace field and quaternion algebra are equal to $K$ and $A$?
\end{question}
This natural question had been discussed for some time, but the complete answer was unknown. Recently J. Kahn and V. Markovic announced a partial answer as follows:

\begin{theorem}\cite{K}\label{KM}
Let $K$ be any real number field and $A$ be any quaternion algebra over $K$ such that $A\otimes_K\mathbb{R}\cong M_2(\mathbb{R})$ and $A\ncong \bigg(\dfrac{1,1}{K}\bigg)$.\footnote{We believe the second condition on $A$ (i.e. $A\ncong \bigg(\dfrac{1,1}{K}\bigg)$) is removable.} Then there exists a closed surface $S_g$ and a hyperbolic structure on it such that $K$ and $A$ are its invariant trace field and quaternion algebra.
\end{theorem}
Their proof uses the techniques developed in their recent proofs of two deep conjectures, the Surface Subgroup Conjecture and the Ehrenpreis Conjecture. 
Although Theorem \ref{KM} answers Question \ref{qe1} partially, an important feature of this theorem is that they realize surfaces via integral traces. 
That is, all the traces in their construction are algebraic integers.\footnote{However, answering the question via integral traces seems to be impossible. See Section \ref{Conjecture}.} 

The aim of this paper is to provide a complete answer to Question \ref{qe1}. In other words, we prove the following theorem:

\begin{theorem}\label{main1}
Let $S_g$ ($g\geq 2$) be any closed surface of genus $g$. Let $K$ be any real number field and $A$ be any quaternion algebra over $K$ such that $A\otimes_K\mathbb{R}\cong M_2(\mathbb{R})$. 
Then there exists a hyperbolic structure on $S_g$ such that $K$ and $A$ arise as its invariant trace field and invariant quaternion algebra.
\end{theorem}
Note that since invariant trace fields and quaternion algebras are invariants of commensurability classes, once we prove the above theorem for $g=2$, the rest of the cases easily follow by looking at covering surfaces. 

We can also ask a similar realization question about the (usual) trace field and quaternion algebra. 
In fact, we prove the following:
\begin{theorem}\label{usual}
Let $S_2$ be a genus $2$ closed surface. Let $K$ be any real number field and $A$ be any quaternion algebra over $K$ such that $A\otimes_K\mathbb{R}\cong M_2(\mathbb{R})$. 
Then there exists a hyperbolic structure on $S_2$ such that $K$ and $A$ arise as its trace field and quaternion algebra. 
Moreover, the invariant trace field and invariant quaternion algebra of it are equal to the trace field and quaternion algebra of it respectively. 
\end{theorem}
We split the proof of this theorem into two parts. We prove the first statement in Section \ref{toward} and, using the proof, we show the second statement in Section \ref{final}. 

Our proofs are based on explicit computations as well as some elementary facts in number theory. 
The basic idea is to create a genus two Riemann surface by attaching two identical copies of a once-punctured torus. 
Then using the work of T. Gauglhofer \cite{TG}, we prove that its trace field is fairly simple, generated only by the traces of three elements. 
We next convert the problem into a system of Diophantine equations and solve these equations. 
The whole process is completely elementary and natural. 

\subsection{Acknowledgments}
I thank Walter Neumann for suggesting this question as well as providing many helpful conversations, and Jeremy Kahn for showing his interest in this work. 
I also thank Ilya Kofman and the anonymous referee for their careful reading and valuable comments. 
Lastly I would like to thank Vlad Markovic and Ian Agol for helpful correspondence and pointing out a mistake in the statement of the main theorem in an earlier version of the paper.  
 
\section{Preliminaries}\label{pre}

In Sections \ref{tf} and \ref{qa}, we quickly review basic definitions and facts which will be used below. For more details, see \cite{RM}. 

\subsection{Trace Field}\label{tf}

Let $\Gamma \subset \text{PSL}(2, \mathbb{R})$ ($\cong\text{Isom}^+(\mathbb{H}^2)$) be a Fuchsian group such that $\mathbb{H}^2/\Gamma$ is a closed hyperbolic surface. 
Let $\bar{\Gamma}\subset \text{SL}(2,\mathbb{R})$ be the inverse image of $\Gamma$ under the projection $\text{SL}(2,\mathbb{R})\rightarrow \text{PSL}(2,\mathbb{R})$. 
Then the \textit{trace field} of $\Gamma$ is defined by
\begin{equation*}
\mathbb{Q}(\{\text{tr}\;\gamma\;|\;\gamma\in\bar{\Gamma}\}).
\end{equation*}
This set is known to be a finite extension number field. For instance, if $\bar{\Gamma}=\;<\gamma_1,\cdots,\gamma_n\;|\;->$, then its trace field is generated by the traces of the following elements:
\begin{equation}\label{ge}
\{\gamma_i,\; \gamma_{j_1}\gamma_{j_2},\; \gamma_{k_1}\gamma_{k_2}\gamma_{k_3}\;|\;
1\leq i\leq n, 1\leq j_1<j_2\leq n, 1\leq k_1<k_2<k_3\leq n\}.
\end{equation}

The \textit{invariant trace field} of $\Gamma$ is defined by 
\begin{equation*}
\mathbb{Q}(\{\text{tr}\;\gamma^2\;|\;\gamma\in\bar{\Gamma}\}),
\end{equation*}
and it is an invariant of the commensurability class of $\Gamma$.

\subsection{Quaternion Algebra}\label{qa}
Let $K$ be a number field. A \textit{quaternion algebra} over $K$ is a four dimensional algebra with basis $\{\bold{1}, \bold{i}, \bold{j}, \bold{k}\}$ such that
\begin{equation*}
\bold{i}^2=a,\; \bold{j}^2=b,\; \bold{i}\bold{j}=-\bold{j}\bold{i}=\bold{k},
\end{equation*}
for some $a,b\in K^*$ and this is often denoted by $\Big(\dfrac{a,\;b}{K}\Big)$. The following equivalence relations are well known: 
\begin{equation}\label{eq}
\Big(\dfrac{a,b}{K}\Big)\cong\Big(\dfrac{au^2,bv^2}{K}\Big)\cong\Big(\dfrac{au^2, bv^2-abw^2}{K}\Big) 
\end{equation}
where $u,v,w\in K$ such that $u,v,bv^2-abw^2\neq 0$.\footnote{To show the equivalence between $A_1=\Big(\dfrac{a,b}{K}\Big)$ and $A_2=\Big(\dfrac{au^2, bv^2-abw^2}{K}\Big)$, let 
$\varphi:A_1\rightarrow A_2$ be a map such that $\varphi(\bold{1})=\bold{1},\varphi(\bold{i})=u\bold{i}$ and $\varphi(\bold{j})=u\bold{j}+w\bold{k}$. Then $\varphi(\bold{j})\varphi(\bold{k})=-\varphi(\bold{k})\varphi(\bold{j})$ and 
so it can be naturally extended to an isomorphism between $A_1$ and $A_2$.}

Let $\Gamma$ be as in Section \ref{tf}. The \textit{quaternion algebra} of $\Gamma$ is defined by
\begin{equation*}
\Big\{\sum_{i=1}^n a_i\gamma_i\;|\; a_i\in \text{Trace field of } \Gamma,\;\;\gamma_i\in \bar{\Gamma} \Big\}.
\end{equation*}
This set is a quaternion algebra over the trace field of $\Gamma$ and equivalent to
\begin{equation}\label{eq.10}
\Big(\dfrac{\text{tr}^2\;\gamma_1-4,\;\; \text{tr}\;[\gamma_1, \gamma_2]-2}{\text{Trace field of }\Gamma} \Big)
\end{equation}
where $\gamma_1, \gamma_2 \in \bar{\Gamma}$ are two hyperbolic elements such that $<\gamma_1,\gamma_2>$ is irreducible. Note that since $\gamma_1$ is hyperbolic, $\text{tr}\;\gamma_1>2$, and thus 
\begin{equation}
\Big(\dfrac{\text{tr}^2\;\gamma_1-4,\;\; \text{tr}\;[\gamma_1, \gamma_2]-2}{K} \Big)\otimes_K\mathbb{R}=
\Big(\dfrac{\text{tr}^2\;\gamma_1-4,\;\; \text{tr}\;[\gamma_1, \gamma_2]-2}{\mathbb{R}} \Big)\cong 
\Big(\dfrac{1,1}{\mathbb{R}} \Big)=M_2(\mathbb{R})
\end{equation}
where $K$ is the trace field of $\Gamma$. The \textit{invariant quaternion algebra} of $\Gamma$ is the algebra generated over the invariant trace field of $\Gamma$ by the squares of the elements of $\bar{\Gamma}$. 
This set is also an invariant of the commensurability class of $\Gamma$.

\subsection{Trace Field Coordinates of Riemann Surfaces}
Each $\gamma\in \text{PSL}(2,\mathbb{R})$ has two inverse images in $\text{SL}(2,\mathbb{R})$. 
We denote the one having the positive trace by $\gamma_+$ and the one having the negative trace by $\gamma_-$.

Following the same notation given in \cite{TG}, we use the matrices below as a basis of $\text{SL}(2,\mathbb{R})$:
\begin{gather*}
\bold{1}=
\begin{pmatrix} 
1 & 0 \\ 
0 & 1 
\end{pmatrix},
\qquad
\bold{I}=
\begin{pmatrix} 
1 & 0 \\ 
0 & -1 
\end{pmatrix},
\qquad
\bold{J}=
\begin{pmatrix} 
0 & 1 \\ 
1 & 0 
\end{pmatrix},
\qquad
\bold{K}=
\begin{pmatrix} 
0 & 1 \\ 
-1 & 0 
\end{pmatrix}.
\end{gather*}
For example, in terms of the above basis, 
$\begin{pmatrix} 
a & b \\ 
c & d 
\end{pmatrix}$
is represented by 
\begin{gather*}
\Big(\dfrac{a+d}{2}\Big)\bold{1}+\Big(\dfrac{a-d}{2}\Big)\bold{I}
+\Big(\dfrac{b+c}{2}\Big)\bold{J}+\Big(\dfrac{b-c}{2}\Big)\bold{K}. 
\end{gather*}

\begin{figure}[htb]
\centering
\includegraphics[width=0.6\linewidth]{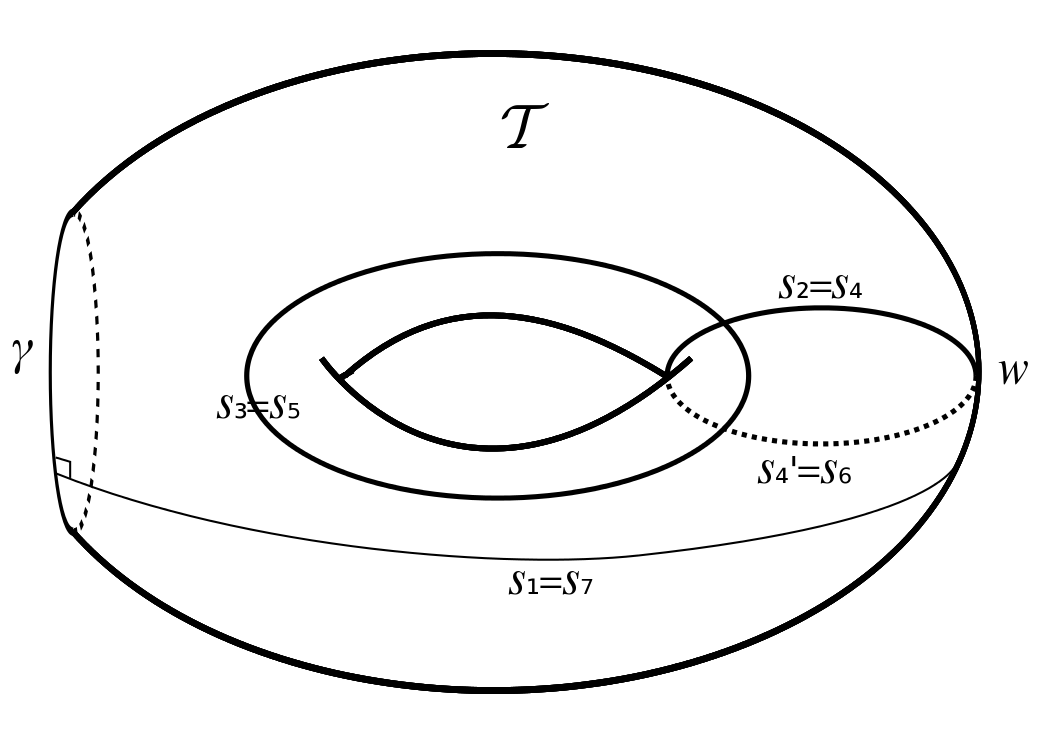}
\caption{}
\label{fig:1}
\end{figure}

\begin{figure}[htb]
\centering
\includegraphics[width=0.8\linewidth]{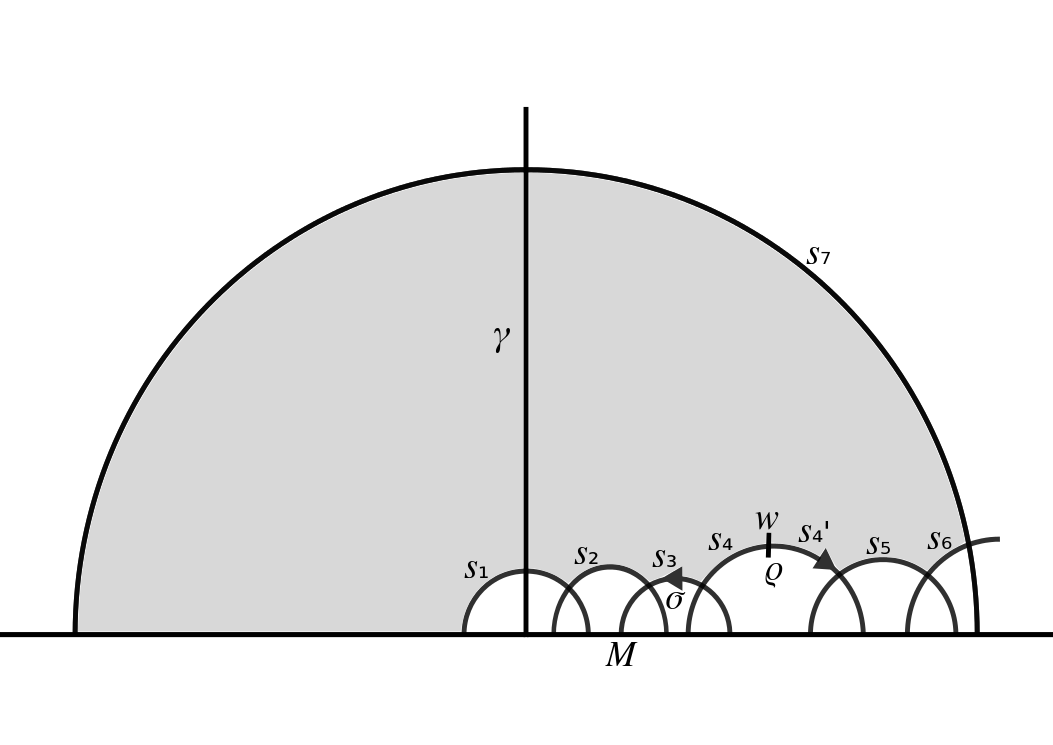}
\caption{$\rho$ and $\sigma$ fix the geodesics containing $s_4\cup s_4'$ and $s_3$ (respectively), and $[\rho,\sigma]$ fixes the $y$-axis. 
The side pairings are given as follows: $\rho(s_3)=s_5, \sigma(s_4)=s_2, \rho\sigma\rho^{-1}(s_4')=s_6$ and $[\rho, \sigma](s_1)=s_7$.}
\label{fig:2}
\end{figure}

Let $\mathcal{T}$ be a hyperbolic once punctured torus whose fundamental domain and picture are shown in Figure \ref{fig:2} and Figure \ref{fig:1} respectively. 
Let $\rho$ and $\sigma$ be the hyperbolic elements given in Figure \ref{fig:2}. Then $\rho_+$ and $\sigma_+$ can be represented in terms of the traces of 
$\rho_+,\sigma_+,(\rho\sigma)_+, [\rho,\sigma]_+$ plus the attracting fixed point of $\sigma$ as follows:\footnote{Note that $[\rho_+,\sigma_+]=-[\rho,\sigma]_+=[\rho,\sigma]_-$.} 

\begin{theorem}\cite{TG}\label{RC} 
Let $\text{tr}\;\rho_+=2r, \text{tr}\;\sigma_+=2s, \text{tr}\;(\rho\sigma)_+=2t,\text{tr}\;([\rho,\sigma]_+)=2c$, 
and $M>0$ be the attracting fixed point of $\sigma$ (as shown in the picture). Then
\allowdisplaybreaks
\begin{equation}\label{ex.1}
\begin{aligned}
&\rho_+=r\bold{1}+\dfrac{r(c+1)}{\sqrt{c^2-1}}\bold{I}-\dfrac{\tilde{M}\big(2rs-t+(c+\sqrt{c^2-1})t\big)}{2(c-1+2s^2)}(\bold{J}+\bold{K})+\dfrac{2rs-t+\big(c-\sqrt{c^2-1}\big)t}{2\tilde{M}(c-1)}\big(\bold{J}-\bold{K}\big),\\
&\sigma_+=s\bold{1}-\dfrac{s(c+1)}{\sqrt{c^2-1}}\bold{I}+\dfrac{\tilde{M}}{2}(\bold{J}+\bold{K})-\dfrac{c-1+2s^2}{2\tilde{M}(c-1)}(\bold{J}-\bold{K}),\\
\end{aligned}
\end{equation}
where
\begin{align}\label{in}
\tilde{M}=M\dfrac{s\sqrt{c^2-1}+(c-1)\sqrt{s^2-1}}{c-1},\quad c=4rst-2r^2-2s^2-2t^2+1\quad\text{and}\quad c,r,s,t>1,M>0.
\end{align}
\end{theorem}
\quad\\
\textbf{Remark.} In \cite{TG}, it is proved that, for any $c,r,s,t,M,\tilde{M}$ satisfying \eqref{in}, the group generated by 
$\rho_+$ and $\sigma_+$ act discretely on $\mathbb{H}$ and generate a fundamental domain as given in Figure \ref{fig:2}. \\

\subsection{Creating a Surface $\mathcal{S}$}
\begin{figure}[htb]
\centering
\includegraphics[width=0.9\linewidth]{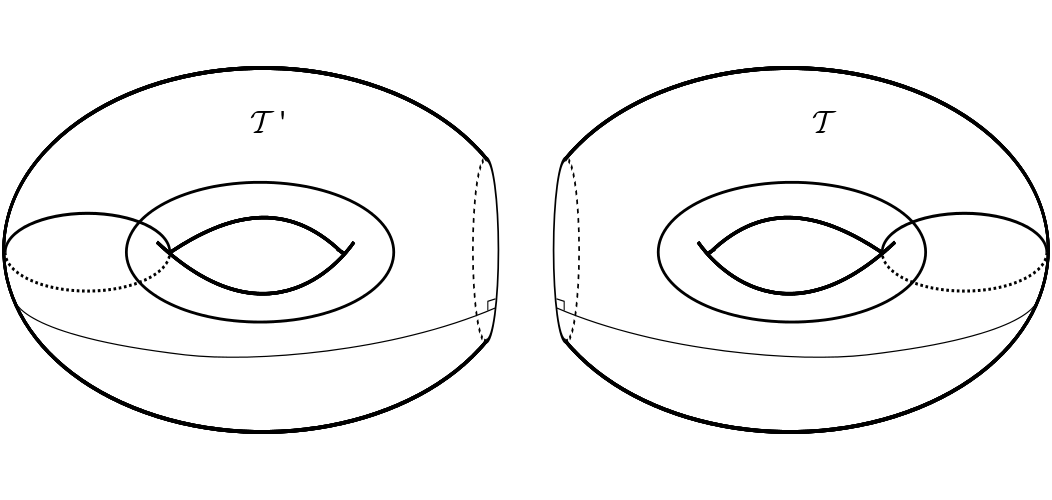}
\caption{}
\label{fig:3}
\end{figure}

\begin{figure}[htb]
\centering
\includegraphics[width=0.8\linewidth]{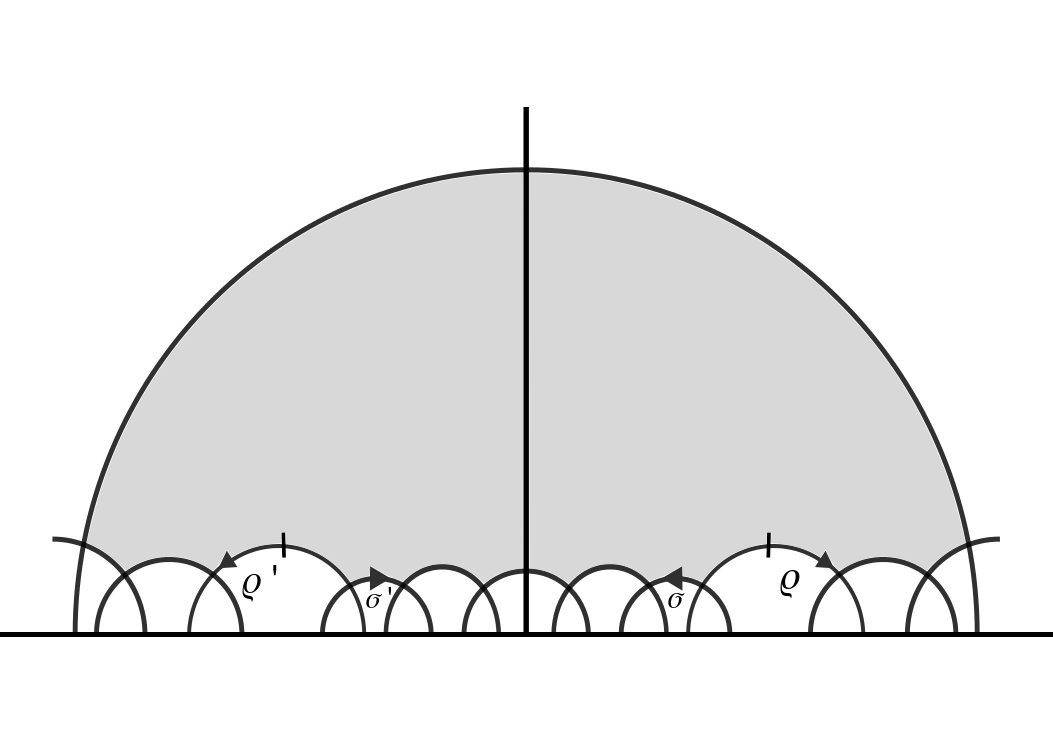}
\caption{}
\label{fig:4}
\end{figure}

Now we create a genus $2$ surface $\mathcal{S}$ by attaching $\mathcal{T}$ to identical symmetric image $\mathcal{T'}$ along their common boundary as shown in Figure \ref{fig:3}. 
(That is, $\mathcal{S}$ is the doubling of $\mathcal{T}$. See Figure \ref{fig:4} for the corresponding fundamental domain.) 
Let $\rho'$ and $\sigma'$ be the elements corresponding to $\rho$ and $\sigma$ (Figure \ref{fig:4}). 
Then $(\rho')_+$ and $(\sigma')_+$ are as follows:\footnote{If $\dfrac{az+b}{cz+d}$ is a M\"{o}bius transformation representing $\rho$, then it is easy to check that 
$\rho'$ is of the form $\dfrac{az-b}{-cz+d}$. In terms of the basis we introduced, 
$\begin{pmatrix} 
a & b \\ 
c & d 
\end{pmatrix}$
is equal to  
$\Big(\dfrac{a+d}{2}\Big)\bold{1}+\Big(\dfrac{a-d}{2}\Big)\bold{I}+\Big(\dfrac{b+c}{2}\Big)\bold{J}+\Big(\dfrac{b-c}{2}\Big)\bold{K}$ 
and $\begin{pmatrix} 
a & -b \\ 
-c & d 
\end{pmatrix}$ is equal to  
$\Big(\dfrac{a+d}{2}\Big)\bold{1}+\Big(\dfrac{a-d}{2}\Big)\bold{I}-\Big(\dfrac{b+c}{2}\Big)\bold{J}-\Big(\dfrac{b-c}{2}\Big)\bold{K}.$ Thus the formula follows.}
\allowdisplaybreaks
\begin{equation}\label{ex.2}
\begin{aligned}
&(\rho')_+=r\bold{1}+\dfrac{r(c+1)}{\sqrt{c^2-1}}\bold{I}+\dfrac{\tilde{M}\big(2rs-t+(c+\sqrt{c^2-1})t\big)}{2(c-1+2s^2)}(\bold{J}+\bold{K})-\dfrac{2rs-t+\big(c-\sqrt{c^2-1}\big)t}{2\tilde{M}(c-1)}\big(\bold{J}-\bold{K}\big),\\
&(\sigma')_+=s\bold{1}-\dfrac{s(c+1)}{\sqrt{c^2-1}}\bold{I}-\dfrac{\tilde{M}}{2}(\bold{J}+\bold{K})+\dfrac{c-1+2s^2}{2\tilde{M}(c-1)}(\bold{J}-\bold{K}),
\end{aligned}
\end{equation}
where
\begin{align*}
\tilde{M}=M\dfrac{s\sqrt{c^2-1}+(c-1)\sqrt{s^2-1}}{c-1},\quad c=4rst-2r^2-2s^2-2t^2+1\quad\text{and}\quad c, r,s,t>1,M>0.
\end{align*}



\subsection{Trace field of $\mathcal{S}$ }


By slightly abusing notation, from now on, we denote $\rho_+, \sigma_+, (\rho')_+$ and $(\sigma')_+$ by $\rho, \sigma, \rho'$ and $\sigma'$ respectively. 
By \eqref{ge}, the trace field of $\mathcal{S}$ is generated by the traces of the following elements:
\begin{equation}\label{eq.1}
\rho,\sigma,\rho',\sigma',\rho\sigma, \rho'\sigma', \rho\rho', \rho\sigma', 
\sigma\rho', \sigma\sigma', 
\rho\sigma\rho', \rho\sigma\sigma', \rho\rho'\sigma', \sigma\rho'\sigma'.
\end{equation}
By \eqref{ex.1} and \eqref{ex.2}, all the elements of $\rho, \sigma, \rho'$ and $\sigma'$ are contained in 
\begin{equation}\label{bo}
\mathbb{Q}(r,s,t, M,\sqrt{s^2-1},\sqrt{c^2-1})\;\;\text{where}\;\;c=4rst-2r^2-2s^2-2t^2+1, 
\end{equation}
and so the elements of \eqref{eq.1} are contained in \eqref{bo} as well. These imply the following simple lemma:
\begin{lemma}
The trace field of $\mathcal{S}$ is contained in \eqref{bo}.
\end{lemma}
Now we further claim the following:
\begin{theorem}\label{thm1}
The trace field of $\mathcal{S}$ is simply equal to $\mathbb{Q}(r,s,t)$.
\end{theorem}
In other words, the trace of any element  in \eqref{eq.1} is represented as a rational function of $r,s, t$ without involving $M, \sqrt{s^2-1}, \sqrt{c^2-1}$. 
Since the proof of this theorem is elementary but lengthy, we postpone it until Section \ref{thms}. 
\section{Proof of Theorem \ref{usual} (Part I)}\label{toward}

In this section, we prove the first part of Theorem \ref{usual}. That is we show that any real number field and any question algebra over it can be realized as the trace field and the quaternion algebra of $\mathcal{S}$. 

By Theorem \ref{thm1}, the trace field of $\mathcal{S}$ is equal to
\begin{equation}\label{tf10}
\mathbb{Q}(r,s,t), 
\end{equation}
and, by \eqref{eq.10}, the quaternion algebra of $\mathcal{S}$ is equivalent to 
\begin{equation}\label{qa10}
\left(\dfrac{(2r)^2-4,\;-2c-2}{\mathbb{Q}(r,s,t)}\right)
\end{equation}
where $c=4rst-2r^2-2s^2-2t^2+1$ and $c,r,s,t>1$.\footnote{Following the same notation given in Section \ref{pre}, recall that $\text{tr}\;([\rho_+,\sigma_+])=-\text{tr}\;([\rho,\sigma]_+)=-2c$.}

Let $K$ be any real number field and $A$ be any quaternion algebra over $K$ such that $A\otimes_K \mathbb{R} \cong M_2(\mathbb{R})$. We suppose $A$ is of the form $\left(\dfrac{a,b}{K}\right)$ for some $a,b\in K$.\footnote{Since  $A\otimes_K \mathbb{R} \cong M_2(\mathbb{R})$, at least one of $a$ and $b$ is positive, 
and, without loss of generality, we assume $a>0$.} 
Since 
\begin{equation*}
\left(\dfrac{a,b}{K}\right)\cong\left(\dfrac{au^2,\;bv^2-abw^2}{K}\right)
\end{equation*}
for any $u,v,w$ (see \eqref{eq}), to achieve the goal, it is enough to find $r,s,t,u,v,w\in \mathbb{R}$ such that
\begin{equation*}
K=\mathbb{Q}(r,s,t)\;\;\text{and}\;\;\left(\dfrac{(2r)^2-4,\;-2c-2}{\mathbb{Q}(r,s,t)}\right)=\left(\dfrac{au^2,\;bv^2-abw^2}{K}\right)
\end{equation*}
where $c=4rst-2r^2-2s^2-2t^2+1$ and $c,r,s,t>1$. In other words, for any given number field $K$ and $a,b\in K$ ($a>0$), it is enough to find $r,s,t,u,v,w\in K$ such that
\begin{equation}\label{eq.11}
\begin{gathered}
c=4rst-2r^2-2s^2-2t^2+1,\\
-2c-2=bv^2-abw^2,\\
(2r)^2-4=au^2\; (a>0),\\
\mathbb{Q}(r,s,t)=K,\\
c,r,s,t>1.
\end{gathered}
\end{equation}
Thus proving the first statement in Theorem \ref{usual} is reduced to solving the above system of Diophantine equations. To simplify the notation, we let $2r=x', 2s=y', 2t=z, 2c=-c'$ and rewrite \eqref{eq.11} as follows:
\begin{gather}
(x')^2+(y')^2+z^2-x'y'z=c'+2,\label{eq.13}\\
c'-2=bv^2-abw^2,\label{eq.15}\\
z^2-4=au^2\; (a>0),\label{eq.14}\\
\mathbb{Q}(x',y',z)=K,\label{eq.100}\\
x',y',z>2, c'<-2.\label{eq.1000}
\end{gather}
To further simplify, we combine \eqref{eq.13} and \eqref{eq.15}, and transform the resulting equation as follows: 
\allowdisplaybreaks
\begin{align*}
&(x')^2+(y')^2+z^2-x'y'z=c'+2 \\
(\text{by}\; \eqref{eq.15})\Rightarrow\; &(x')^2+(y')^2+z^2-x'y'z=bv^2-abw^2+4\\
\Rightarrow\; &(x'-\dfrac{y'z}{2})^2-\dfrac{(y')^2z^2}{4}+(y')^2+z^2=bv^2-abw^2+4\\
\Rightarrow\; &(x'-\dfrac{y'z}{2})^2-\dfrac{(y')^2}{4}(z^2-4)+z^2=bv^2-abw^2+4\\
(\text{by}\; \eqref{eq.14})\Rightarrow\; &(x'-\dfrac{y'z}{2})^2-\dfrac{(y')^2}{4}(au^2)+au^2+4=bv^2-abw^2+4\\
\Rightarrow\; &(x'-\dfrac{y'z}{2})^2-\dfrac{(y')^2}{4}(au^2)+au^2=bv^2-abw^2\\
\Rightarrow\; &(x'-\dfrac{y'z}{2})^2-\dfrac{(y')^2}{4}(au^2)+au^2-bv^2+abw^2=0.
\end{align*}
Let $x=x'-\dfrac{y'z}{2}$ and $y=\dfrac{y'u}{2}$. Then \eqref{eq.13} - \eqref{eq.100} are reduced to 
\begin{equation}\label{lll}
\begin{gathered}
x^2-ay^2+au^2-bv^2+abw^2=0,\\
z^2-4=au^2\; (a>0),\\
\mathbb{Q}(x,\dfrac{y}{u},z)=K\;(u\neq 0),
\end{gathered}
\end{equation}
and \eqref{eq.1000} is equivalent to 
\begin{gather}\label{eq.1010}
x+\dfrac{yz}{2u}>2,\;\dfrac{y}{u}>1,\;z>2, \;abw^2-bv^2>4.
\end{gather}
Now the first statement in Theorem \ref{usual} follows from the theorem below: 
\begin{theorem}\label{main3}
For any real number field $K$ and any $a,b\in K$ $(a>0)$, there exists 
$x,y,z,u,v,w\in K$ satisfying the following:
\begin{gather}
x^2-ay^2+au^2-bv^2+abw^2=0,\label{eq.16}\\
z^2-4=au^2\; (a>0),\label{eq.20}\\
\mathbb{Q}\Big(x,\dfrac{y}{u},z\Big)=K,\label{eq.29}\\
u\neq0,\;x+\dfrac{yz}{2u}>2,\;\dfrac{y}{u}>1,\;z>2, \;abw^2-bv^2>4.\label{eq.1111}
\end{gather}
\end{theorem}
Here is the outline of the proof. We resolve each equation step by step. 
First, in \textbf{Step 1}, we parameterize solutions of \eqref{eq.16}, using five variables $m_i$ ($1\leq i\leq 5$), as given in \eqref{eq.19}. 
Since \eqref{eq.16} is homogeneous equation of degree $2$, any constant multiple of \eqref{eq.19} is again a solution of \eqref{eq.16}, and thus the complete set of solutions is of the form in \eqref{100} with six degrees of freedom. 

Next, in \textbf{Step 2}, we resolve the second equation \eqref{eq.20}. Among those six variables obtained in the previous step, by fixing one as a function (this function is induced from \eqref{eq.20}) of the rest five variables, $m_i$($1\leq i\leq 5$), 
we let \eqref{eq.20} always hold for any $m_i$. In other words, instead of losing one variable, we remove \eqref{eq.20} in return.  

In \textbf{Step 3}, we deal with \eqref{eq.29}. This is the most delicate one. Intuitively, for any given real number field $K$, the set of generators of $K$ is dense in $\mathbb{R}$. 
But, in our case, since each of $x,\dfrac{y}{u}$ and $z$ is a rational function of $m_i$, we need to show there is a particular way to choose $m_i$ so that one of $x,\dfrac{y}{u}$ and $z$ generates $K$. 
Specifically, in Lemma \ref{key0}, we use $\dfrac{y}{u}$ as a generator of $K$ and show 
\begin{equation*}
 \mathbb{Q}\Big(x,\dfrac{y}{u},z\Big)=\mathbb{Q}\Big(\dfrac{y}{u}\Big)=K.
\end{equation*}
In fact, the lemma will be further reduced to a simpler but equivalent statement (i.e. Lemma \ref{key}), and this will be proved using subsequent claims. 

In \textbf{Step 4}, it is shown one can further control the ranges of $m_i$ to get a solution satisfying \eqref{eq.1111} as well.

\begin{proof}
\textbf{Step 1.} Having $(x,y,u,v,w)=(0,1,1,0,0)$ is as an initial solution, we employ a well-known technique in number theory to find other solutions of \eqref{eq.16}. 
Let 
\begin{equation}
\begin{aligned}\label{eq.18}
&\dfrac{x}{m_1}=\dfrac{y-1}{m_2}=\dfrac{u-1}{m_3}=\dfrac{v}{m_4}=\dfrac{w}{m_5}=t.
\end{aligned}
\end{equation}
Combining \eqref{eq.18} with \eqref{eq.16}, we get 
\allowdisplaybreaks
\begin{align}
&(m_1t)^2-a(m_2t+1)^2+a(m_3t+1)^2-b(m_4t)^2+ab(m_5t)^2=0\notag\\
\Rightarrow\; &(m_1t)^2-a\big((m_2t)^2+2m_2t+1\big)+a\big((m_3t)^2+2m_3t+1\big)-b(m_4t)^2+ab(m_5t)^2=0\notag\\
\Rightarrow\; &(m_1t)^2-a\big((m_2t)^2+2m_2t\big)+a\big((m_3t)^2+2m_3t\big)-b(m_4t)^2+ab(m_5t)^2=0\notag\\
\Rightarrow\; &\big(m_1^2-am_2^2+am_3^2-bm_4^2+abm_5^2\big)t^2-2a\big(m_2-m_3\big)t=0\notag\\
\Rightarrow\; &\big(m_1^2-am_2^2+am_3^2-bm_4^2+abm_5^2\big)t=2a\big(m_2-m_3\big)\notag\\
\Rightarrow\; &t=\dfrac{2a\big(m_2-m_3\big)}{m_1^2-am_2^2+am_3^2-bm_4^2+abm_5^2}.\label{eq.11111}
\end{align}
So, for any $m_1, m_2, m_3, m_4, m_5\in K$ such that 
\begin{equation}\label{eq.0101}
m_1^2-am_2^2+am_3^2-bm_4^2+abm_5^2\neq 0,
\end{equation}
the following is a solution of \eqref{eq.16}: 
\allowdisplaybreaks
\begin{equation}\label{eq.19}
\begin{aligned}
&x=m_1t, & &y=\left(m_2t+1\right), & &u=\left(m_3t+1\right), & &v=m_4t, & &w=m_5t
\end{aligned}
\end{equation}
where $t$ is the one given in \eqref{eq.11111}. 
Furthermore since \eqref{eq.16} is homogeneous, for any $d\in K$, the following is also a solution of \eqref{eq.16}:
\allowdisplaybreaks
\begin{equation}\label{100}
\begin{aligned}
&x=m_1td,& &y=\left(m_2t+1\right)d, & &u=\left(m_3t+1\right)d, & &v=m_4td, & &w=m_5td,
\end{aligned}
\end{equation}
where $m_i$ ($1\leq i\leq 5$) and $t$ are the same as above. 
\\
\\
\textbf{Step 2.} Using the common factor $d$ in \eqref{100}, we now resolve the second equation \eqref{eq.20}. 
Since 
\begin{equation*}
z^2-4=au^2=a(m_3t+1)^2d^2,
\end{equation*}
by letting
\begin{gather*}
z+2=a(m_3t+1)d\quad \text{and}\quad z-2=(m_3t+1)d,
\end{gather*}
we get 
\begin{gather*}
(z+2)-(z-2)=4=d(m_3t+1)(a-1)\\
\end{gather*}
and so
\begin{gather}\label{co}
d=\dfrac{4}{(m_3t+1)(a-1)}.
\end{gather}
Thus, by fixing $d$ as in \eqref{co}, we can always make \eqref{eq.20} hold for any $m_i$ ($1\leq i\leq 5$) satisfying \eqref{eq.0101} and $m_3t+1\neq 0$.\footnote{By \eqref{eq}, we can assume $a>1$. We will talk more about this in \textbf{Step 4}.}  
\\
\\
\textbf{Step 3.} Next we resolve the third condition \eqref{eq.29}. Among three elements $x,\dfrac{y}{u}, z$, 
we use $\dfrac{y}{u}$ as a generator of $K$. That is, we show the following lemma:
\begin{lemma}\label{key0}
Let $x,y,z,u,v$ be as given in \eqref{100}. For any real number field $K$, there exists $m_1,m_2,m_3,m_4,m_5\in K$ such that 
\begin{equation}
\mathbb{Q}\Big(x,\dfrac{y}{u},z\Big)=\mathbb{Q}\Big(\dfrac{y}{u}\Big)=K.
\end{equation}
More precisely, there exists a generator $g$ of $K$ and $m_1,m_2,m_3,m_4,m_5\in K$ such that
\begin{equation}\label{f}
\dfrac{y}{u}=\dfrac{m_2t+1}{m_3t+1}=g,
\end{equation}
where $t=\dfrac{2a(m_2-m_3)}{m_1^2-am_2^2+am_3^2-bm_4^2+abm_5^2}$. 
\end{lemma}
Before proving the lemma, we first simplify \eqref{f} as follows:
\allowdisplaybreaks
\begin{align}
&\dfrac{m_2t+1}{m_3t+1}=g\notag\\
\Rightarrow\; &m_2t+1=g(m_3t+1)\notag\\
\Rightarrow\; &\dfrac{2a(m_2-m_3)m_2}{m_1^2-am_2^2+am_3^2-bm_4^2+abm_5^2}+1=g\bigg(\dfrac{2a(m_2-m_3)m_3}{m_1^2-am_2^2+am_3^2-bm_4^2+abm_5^2}+1\bigg)\notag\\
\Rightarrow\; &\dfrac{2a(m_2-m_3)m_2+\big(m_1^2-am_2^2+am_3^2-bm_4^2+abm_5^2\big)}{m_1^2-am_2^2+am_3^2-bm_4^2+abm_5^2}\notag\\
&=g\bigg(\dfrac{2a(m_2-m_3)m_3+\big(m_1^2-am_2^2+am_3^2-bm_4^2+abm_5^2\big)}{m_1^2-am_2^2+am_3^2-bm_4^2+abm_5^2}\bigg)\notag\\
\Rightarrow\; &2a(m_2-m_3)m_2+\big(m_1^2-am_2^2+am_3^2-bm_4^2+abm_5^2\big)\notag\\
&=g\Big(2a(m_2-m_3)m_3+\big(m_1^2-am_2^2+am_3^2-bm_4^2+abm_5^2\big)\Big)\notag\\
\Rightarrow\; &(g-1)\big(m_1^2-am_2^2+am_3^2-bm_4^2+abm_5^2\big)+2ga(m_2m_3-m_3^2)-2a(m_2^2-m_2m_3)=0\notag\\
\Rightarrow\; &(g-1)m_1^2-(g-1)bm_4^2+(g-1)abm_5^2-(g-1)am_2^2+(g-1)am_3^2+2gam_2m_3-2gam_3^2\notag\\
&-2am_2^2+2am_2m_3=0\notag\\
\Rightarrow\; &(g-1)m_1^2-(g-1)bm_4^2+(g-1)abm_5^2-(g+1)am_2^2-(g+1)am_3^2+2gam_2m_3+2am_2m_3=0\notag\\
\Rightarrow\; &(g-1)m_1^2-(g-1)bm_4^2+(g-1)abm_5^2-(g+1)a(m_2^2-2am_2m_3+m_3^2)=0\notag\\
\Rightarrow\; &(g-1)m_1^2-(g-1)bm_4^2+(g-1)abm_5^2-(g+1)a(m_2-m_3)^2=0\notag\\
\Rightarrow\; &m_1^2-bm_4^2+abm_5^2-\dfrac{a(g+1)}{g-1}(m_2-m_3)^2=0\notag\\
\Rightarrow\; &m_1^2-bm_4^2+abm_5^2-a\Big(1+\dfrac{2}{g-1}\Big)(m_2-m_3)^2=0\notag\\
\Rightarrow\; &\dfrac{1}{a}m_1^2-\dfrac{b}{a}m_4^2+bm_5^2-\Big(1+\dfrac{2}{g-1}\Big)(m_2-m_3)^2=0\notag\\
\Rightarrow\; &\dfrac{1}{a}\Big(\dfrac{m_1}{m_2-m_3}\Big)^2-\dfrac{b}{a}\Big(\dfrac{m_4}{m_2-m_3}\Big)^2+b\Big(\dfrac{m_5}{m_2-m_3}\Big)^2-\Big(1+\dfrac{2}{g-1}\Big)=0.\notag\\
\Rightarrow\; &\dfrac{1}{a}\Big(\dfrac{m_1}{m_2-m_3}\Big)^2-\dfrac{b}{a}\Big(\dfrac{m_4}{m_2-m_3}\Big)^2+b\Big(\dfrac{m_5}{m_2-m_3}\Big)^2=1+\dfrac{2}{g-1}.\label{f1}
\end{align}
Thus the following lemma implies Lemma \ref{key0}:
\begin{lemma}\label{key100}
Let $K$ be any real number field and $a,b$ be arbitrary nonzero elements in $K$ such that $a>0$. Then there exist a generator $g$ of $K$ and $m_1,m_2,m_3,m_4,m_5\in K$ satisfying 
\begin{equation}
\dfrac{1}{a}\Big(\dfrac{m_1}{m_2-m_3}\Big)^2-\dfrac{b}{a}\Big(\dfrac{m_4}{m_2-m_3}\Big)^2+b\Big(\dfrac{m_5}{m_2-m_3}\Big)^2=1+\dfrac{2}{g-1}.\label{f10}
\end{equation}
\end{lemma}
By changing variables as 
\begin{equation}
a'=\dfrac{1}{a},\;\;\;b'=b,\;\;\;x=\dfrac{m_1}{m_2-m_3},\;\;\;z=\dfrac{m_4}{m_2-m_3},\;\;\;y=\dfrac{m_5}{m_2-m_3},\;\;\;g'=1+\dfrac{2}{g-1},
\end{equation}
Lemma \ref{key100} now follows by the lemma below: 
\begin{lemma}\label{key}
Let $K$ be any real number field and $a',b'$ be arbitrary nonzero elements in $K$. Then there exist a generator $g'$ of $K$ and $x,y,z\in K$ satisfying 
\begin{equation}\label{eqi}
a'x^2+b'y^2-a'b'z^2=g'.
\end{equation}
\end{lemma}
Before proving Lemma \ref{key}, we first state two claims which will be used in the proof. 
\begin{claim}\label{lm1}
Let $\alpha\in\mathbb{\bar{Q}}$. If a minimal polynomial of $\alpha$ contains a nonzero term of odd degree, then $\mathbb{Q}(\alpha)=\mathbb{Q}(\alpha^2)$. Thus, for any real number field $K$, the set of generators of $K$ satisfying this property is dense in $\mathbb{R}$. 
\end{claim}
\begin{proof}
Let $f(x)=0$ be a minimal polynomial of $\alpha$. Since it contains a nonzero term of odd degree, we can express $f$ as follows:
\begin{equation*}
f(x)=x(a_{2n}x^{2n}+a_{2n-1}x^{2n-1}+\cdots+a_2x^2+a_0)+(b_{2m}x^{2m}+b_{2m-1}x^{2m-1}+\cdots+b_2x^2+b_0).
\end{equation*}
Since $f(x)$ is a minimal polynomial of $\alpha$, $a_{2n}\alpha^{2n}+a_{2n-1}\alpha^{2n-1}+\cdots+a_2\alpha^2+a_0\neq 0$, and so 
\begin{equation*}
\alpha=-\dfrac{b_{2m}\alpha^{2m}+b_{2m-1}\alpha^{2m-1}+\cdots+b_2\alpha^2+b_0}{a_{2n}\alpha^{2n}+a_{2n-1}\alpha^{2n-1}+\cdots+a_2\alpha^2+a_0},
\end{equation*}
which implies $\mathbb{Q}(\alpha)=\mathbb{Q}(\alpha^2)$. 

Let $\alpha$ be a generator of $K$. Then, for any nonzero $r\in\mathbb{Q}$, either a minimal polynomial of $\alpha$ or a minimal polynomial of $\alpha+r$ contains a nonzero term of odd degree. The second statement easily follows.  
\end{proof}
\begin{claim}\label{lm2}
Let $\alpha,\beta$ be algebraic numbers. Then, except for finitely many $r\in\mathbb{Q}$, it always satisfies $\mathbb{Q}(\alpha+r\beta)=\mathbb{Q}(\alpha,\beta)$. 
\end{claim}
\begin{proof}
This is a well-known fact and used to prove the existence of a primitive element in a number field. See \cite{hun}, for instance, for a proof.  
\end{proof}

\begin{proof}[Proof of Lemma \ref{key}.]
By Claim \ref{lm2}, there exists $(x_0, y_0)\in \mathbb{Q}^2$ such that $\mathbb{Q}(a'x_0^2+b'y_0^2)=\mathbb{Q}(a', b')$. 
Let $F=\mathbb{Q}(a',b')$, and $z_0$ be an element of $K$ such that $K=F(z_0^2)$ obtained from Claim \ref{lm1}. By Claim \ref{lm2}, we can further assume that $x_0$ and $y_0$ satisfy $\mathbb{Q}(a'x_0^2+b'y_0^2-a'b'z_0^2)=\mathbb{Q}(a'x_0^2+b'y_0^2, a'b'z_0^2)$. 
Now the following equalities hold:
\begin{gather*}
\mathbb{Q}(a'x_0^2+b'y_0^2-a'b'z_0^2)=\mathbb{Q}(a'x_0^2+b'y_0^2, a'b'z_0^2)=\mathbb{Q}(a',b')(a'b'z_0^2)=\mathbb{Q}(a',b')(z_0^2)=F(z_0^2)=K.
\end{gather*}
Taking $g'=a'x_0^2+b'y_0^2-a'b'z_0^2$, we get the desired result. 
\end{proof}
\quad\\
\textbf{Step 4.} Now we find a solution satisfying \eqref{eq.1111}. 
For sufficiently small number $\epsilon>0$, using \eqref{eq}, we first assume that $a$ and $b$ satisfy 
\begin{equation}
1<a,\;b<1+\epsilon. 
\end{equation}
Since  $a>1$, we have $z^2>4$ by \eqref{eq.20}. Also if $\dfrac{y}{u}>1$ and $z>2$, then $\dfrac{yz}{2u}>1$. 
Thus \eqref{eq.1111} can be replaced by  
\begin{equation}\label{co2}
u \neq 0,\;\;\; x>1,\;\;\;\dfrac{y}{u}>1,\;\;\;abw^2-bv^2>4. 
\end{equation}
Following the formulas given in \eqref{100}, \eqref{co2} is equivalent to
\begin{equation}\label{co1}
(m_3t+1)d\neq 0,\;\;\;m_1td>1,\;\;\;\dfrac{m_2t+1}{m_3t+1}>1,\;\;\;ab(m_5td)^2-b(m_4td)^2>4.
\end{equation}

By the proof of Lemma \ref{key}, the set of 
\begin{equation*}
\Big(\dfrac{m_1}{m_2-m_3},\dfrac{m_4}{m_2-m_3}, \dfrac{m_5}{m_2-m_3}\Big)  
\end{equation*}
such that 
\begin{equation*}
\dfrac{1}{a}\Big(\dfrac{m_1}{m_2-m_3}\Big)^2-\dfrac{b}{a}\Big(\dfrac{m_4}{m_2-m_3}\Big)^2+b\Big(\dfrac{m_5}{m_2-m_3}\Big)^2 
\end{equation*}
is a generator of $K$ is dense in $\mathbb{R}^3$. Thus for sufficiently large $L>0$, we can pick $m_i$ satisfying
\begin{equation}
L<m_1,\;m_2,\;m_4<L+\epsilon,\;\;\;2L<m_5<2L+\epsilon,\;\;\;0<m_3<\epsilon. 
\end{equation}
Then
\begin{align*}
m_1^2-a(m_2^2-m_3^3)-bm_4^2+abm_5^2&<(L+\epsilon)^2-(L^2-\epsilon^2)-L^2+(1+\epsilon)^2(2L+\epsilon)^2<4L^2\\
\end{align*}
and 
\begin{align*}
m_1^2-a(m_2^2-m_3^3)-bm_4^2+abm_5^2>L^2-(1+\epsilon)(L+\epsilon)^2-(1+\epsilon)(L+\epsilon)^2+(2L)^2>2L^2.
\end{align*}
Thus 
\begin{equation}\label{1}
\dfrac{1}{3L}<\dfrac{1}{2L}\Big(1-\dfrac{\epsilon}{L}\Big)=\dfrac{2(L-\epsilon)}{4L^2}<t=\dfrac{2a(m_2-m_3)}{m_1^2-a(m_2^2-m_3^2)-bm_4^2+abm_5^2}<\dfrac{2(1+\epsilon)(L+\epsilon)}{2L^2}<\dfrac{2}{L}
\end{equation}
and
\begin{equation}\label{2}
1<m_3t+1<\epsilon\dfrac{2}{L}+1.
\end{equation}
Now \eqref{1} and \eqref{2} imply
\begin{equation*}
d=\dfrac{1}{(m_3t+1)(a-1)}>\dfrac{1}{(\frac{2\epsilon}{L}+1)\epsilon}>\dfrac{1}{2\epsilon}.
\end{equation*}
Clearly we have $(m_3t+1)d\neq0$, and
\begin{equation*}
m_1td>L\Big(\dfrac{1}{3L}\Big)\Big(\dfrac{1}{2\epsilon}\Big)=\dfrac{1}{6\epsilon}>1
\end{equation*}
for sufficiently small $\epsilon$. Since $m_2>m_3$, $\dfrac{m_2t+1}{m_3t+1}>1$. Lastly 
\begin{equation*}
ab(m_5td)^2-b(m_4td)^2=t^2d^2(abm_5^2-bm_4^2)>\Big(\dfrac{1}{3L}\Big)^2\Big(\dfrac{1}{2\epsilon}\Big)^2\big(4L^2-(1+\epsilon)(L+\epsilon)^2\big)>\dfrac{1}{18\epsilon^2}>4
\end{equation*}
since $\epsilon$ is sufficiently small.

Thus it satisfies \eqref{co1} and this completes the proof of Theorem \ref{main3}.
\end{proof}
\quad\\
\textbf{Remark.} 
For later use, we introduce a simple criteria to generate infinitely many different solutions from an initial solution. This method will be used in the next section. 
For any real number field $K$ and $a,b\in K$ $(a>0)$, let $m_i=m_{i_{0}}\in K$ ($1\leq i\leq 5$) be numbers obtained by following the above proof of Theorem \ref{main3}. 
(That is, we get  $x,y,z,u,v,w$ satisfying \eqref{eq.16} - \eqref{eq.1111} by letting $m_i=m_{i_{0}}$ in \eqref{100}.) Then 
\begin{equation*}
\dfrac{1}{a}\Big(\dfrac{m_{1_{0}}}{m_{2_{0}}-m_{3_{0}}}\Big)^2-\dfrac{b}{a}\Big(\dfrac{m_{4_{0}}}{m_{2_{0}}-m_{3_{0}}}\Big)^2+b\Big(\dfrac{m_{5_{0}}}{m_{2_{0}}-m_{3_{0}}}\Big)^2
\end{equation*}
is a generator of $K$, and, by \eqref{f1}, we have 
\begin{equation}\label{kkk}
\dfrac{1}{a}\Big(\dfrac{m_{1_0}}{m_{2_0}-m_{3_0}}\Big)^2-\dfrac{b}{a}\Big(\dfrac{m_{4_0}}{m_{2_0}-m_{3_0}}\Big)^2+b\Big(\dfrac{m_{5_0}}{m_{2_0}-m_{3_0}}\Big)^2=1+\dfrac{2}{\frac{y_{0}}{u_0}-1}
\end{equation}
where $y_{0}$ and $u_0$ are the values of $y$ and $u$ when $m_i=m_{i_0}$. For $r\in \mathbb{Q}\backslash \{0\}$, since
\begin{equation*}
\dfrac{1}{a}\Big(\dfrac{rm_{1_{0}}}{m_{2_{0}}-m_{3_{0}}}\Big)^2-\dfrac{b}{a}\Big(\dfrac{rm_{4_{0}}}{m_{2_{0}}-m_{3_{0}}}\Big)^2+b\Big(\dfrac{rm_{5_{0}}}{m_{2_{0}}-m_{3_{0}}}\Big)^2
\end{equation*}
is also a generator of $K$, we get another $x,y,z,u,v,w$ satisfying \eqref{eq.16} - \eqref{eq.29} from
\begin{equation}\label{inf}
(m_1, m_2, m_3, m_4, m_5)=(rm_{1_0},m_{2_0},m_{3_0},rm_{4_0},rm_{5_0}).
\end{equation}
Furthermore, if $r$ is sufficiently close to $1$, then the solution satisfies \eqref{eq.1111} as well. Thus, through this way, we can produce infinitely many different solutions, which are all close to the initial solution 
Finally remark that if $y_r$ and $u_r$ are the values of $y$ and $u$ when $m_i$ ($1\leq i\leq 5$) are as given in \eqref{inf}, then the following two equalities hold:
\begin{equation}\label{kkkk}
\dfrac{1}{a}\Big(\dfrac{rm_{1_0}}{m_{2_0}-m_{3_0}}\Big)^2-\dfrac{b}{a}\Big(\dfrac{rm_{4_0}}{m_{2_0}-m_{3_0}}\Big)^2+b\Big(\dfrac{rm_{5_0}}{m_{2_0}-m_{3_0}}\Big)^2=1+\dfrac{2}{\frac{y_{r}}{u_r}-1},
\end{equation}
\begin{equation}\label{re}
r^2\Big(1+\dfrac{2}{\frac{y_0}{u_0}-1}\Big)=1+\dfrac{2}{\frac{y_r}{u_r}-1}
\end{equation}
(the first one follows from \eqref{f1}, and the second one from \eqref{kkk} and \eqref{kkkk}).  

\section{Proof of Theorem \ref{usual} (Part II)}\label{final}
Now we complete the proof of Theorem \ref{usual}. 
We will show that both the trace field and the invariant trace field of $\mathcal{S}$ can be made equal in the proof of Theorem \ref{main3}. 
Remark that once the trace field is equal to the invariant trace field, then the quaternion algebra is also equal to the invariant quaternion algebra.\footnote{Since both the quaternion algebra and the invariant quaternion algebra are then defined over the same field
and the invariant quaternion algebra is a subalgebra of the quaternion algebra, both are the same.} 

In the previous section, we used the trace of $\sigma$ (recall $\text{tr}\;\sigma=2y/u$) to realize an arbitrary given real number field $K$ as the trace field of $\mathcal{S}$. 
Since the invariant trace field of $\mathcal{S}$ contains $\text{tr}\;\sigma^2$ and
\begin{equation*}
\text{tr}\;\sigma^2=(\text{tr}\;\sigma)^2-4,
\end{equation*}
the following strengthened version of Theorem \ref{main3} implies the last statement of Theorem \ref{usual}:
\begin{theorem}\label{main4}
Theorem \ref{main3} still holds after replacing \eqref{eq.29} by
\begin{equation*}
\mathbb{Q}\Big(\dfrac{y}{u}\Big)=\mathbb{Q}\Big(\dfrac{y^2}{u^2}\Big)=K.
\end{equation*}
\end{theorem}
Intuitively this theorem is clear since we have a dense set of solutions satisfying \eqref{eq.16} - \eqref{eq.1111} in Theorem \ref{main3}. 
But, for the sake of completeness, we will  give the complete proof below. Before proving the theorem, we prove the following lemma first: 
\begin{lemma}\label{key2}
Let $\alpha$ be an algebraic number and $f(x)=a_{n}x^{n}+\cdots+a_0$ be a minimal polynomial of $\alpha$. 
Then $(x-1)^nf\Big(1+\dfrac{2}{x-1}\Big)$ is a minimal polynomial of $1+\dfrac{2}{\alpha-1}$, and the coefficient of $x^{n-1}$ in $(x-1)^nf\Big(1+\dfrac{2}{x-1}\Big)$ is 
\begin{gather*}
\sum^{n}_{i=0} \big(-n+2i\big)a_i.
\end{gather*}
\end{lemma}
\begin{proof}
Let $\beta=1+\dfrac{2}{\alpha-1}$. Then $\alpha=1+\dfrac{2}{\beta-1}$, and so the first part follows easily. 

We next expand $f\Big(1+\dfrac{2}{x-1}\Big)$ as below:
\allowdisplaybreaks
\begin{equation}\label{eq.21}
\begin{aligned}
&a_{n}\Big(1+\dfrac{2}{x-1}\Big)^n+a_{n-1}\Big(1+\dfrac{2}{x-1}\Big)^{n-1}+\cdots+a_1\Big(1+\dfrac{2}{x-1}\Big)+a_0\\
=&a_{n}\bigg(1+n\Big(\dfrac{2}{x-1}\Big)+\cdots+n\Big(\dfrac{2}{x-1}\Big)^{n-1}+\Big(\dfrac{2}{x-1}\Big)^{n}\bigg)\\
&+a_{n-1}\bigg(1+(n-1)\Big(\dfrac{2}{x-1}\Big)+\cdots+(n-1)\Big(\dfrac{2}{x-1}\Big)^{n-2}+\Big(\dfrac{2}{x-1}\Big)^{n-1}\bigg)\\
&\hdots\\
&+a_1\bigg(1+\Big(\dfrac{2}{x-1}\Big)\bigg)+a_0.
\end{aligned}
\end{equation}
Multiplying \eqref{eq.21} by $(x-1)^{n}$, $(x-1)^nf\Big(1+\dfrac{2}{x-1}\Big)$ is equal to    
\allowdisplaybreaks
\begin{equation}\label{su}
\begin{aligned}
&a_{n}\Big((x-1)^{n}+n\cdot 2 (x-1)^{n-1}+\cdots+n\cdot 2^{n-1}(x-1)+2^{n}\Big)\\
&+a_{n-1}\Big((x-1)^{n}+(n-1)2 (x-1)^{n-1}+\cdots+(n-1)2^{n-2}(x-1)^2+2^{n-1}(x-1)\Big)\\
&\hdots\\
&+a_1\Big((x-1)^{n}+2 (x-1)^{n-1}\Big)+a_0(x-1)^{n}.
\end{aligned}
\end{equation}
The coefficient of $x^{n-1}$ in \eqref{su} is 
\begin{gather*}
a_{n}\big(-n+n\cdot 2\big)+a_{n-1}\big(-n+(n-1)2\big)+\cdots+a_1 \big(-n+2\big)+a_0 \big(-n\big),
\end{gather*}
which is equal to
\begin{gather}\label{coe}
\sum^{n}_{i=0} (-n+2i)a_i.
\end{gather}
\end{proof}
Now we are ready to prove Theorem \ref{main4}.
\begin{proof}[Proof of Theorem \ref{main4}.]
Let $y=y_{0}$, $u=u_0$ be the pair given in the remark in the previous section (after the proof of Theorem \ref{main3}). Let $f(x)$ be a minimal polynomial of $1+\dfrac{2}{\frac{y_0}{u_0}-1}$. Then, by Lemma \ref{key2},    
\begin{equation}\label{polyn}
(x-1)^{\text{deg}\;f}f\Big(1+\dfrac{2}{x-1}\Big)
\end{equation}
is a minimal polynomial of $\dfrac{y_0}{u_0}$. If \eqref{polyn} contains a nonzero term of odd degree, then $\mathbb{Q}\Big(\dfrac{y_0}{u_0}\Big)=\mathbb{Q}\Big(\dfrac{y^2_0}{u^2_0}\Big)$ by Claim \ref{lm1}, and so we are done.

Otherwise, suppose the degrees of all the nonzero terms of \eqref{polyn} are even. Let $y_r$ and $u_r$ be the same ones given in the remark in Section \ref{toward}, and $f(x)$ be of the form $a_{2n}x^{2n}+\cdots+a_0$. Then   
\begin{equation}\label{poly}
\begin{aligned}
f\Big(\dfrac{x}{r^2}\Big)=&a_{2n}\Big(\dfrac{x}{r^2}\Big)^{2n}+a_{2n-1}\Big(\dfrac{x}{r^2}\Big)^{2n-1}+\cdots+a_1\Big(\dfrac{x}{r^2}\Big)+a_0\\
=&\Big(\dfrac{a_{2n}}{r^{4n}}\Big)x^{2n}+\Big(\dfrac{a_{2n-1}}{r^{4n-2}}\Big)x^{2n-1}+\cdots+\Big(\dfrac{a_1}{r^2}\Big)x+a_0\\
\end{aligned}
\end{equation}
is a minimal polynomial of $r^2 \Big(1+\dfrac{2}{\frac{y_0}{u_0}-1}\Big)$.
To simplify the notation, we denote the polynomial in \eqref{poly} by $g(x)$ and $r^2 \Big(1+\dfrac{2}{\frac{y_0}{u_0}-1}\Big)$ by $\alpha$. 
Then $\alpha=1+\dfrac{2}{\frac{y_r}{u_r}-1}$ by \eqref{re}, and so $\dfrac{y_r}{u_r}=1+\dfrac{2}{\alpha-1}$. Since $g(x)$ is a minimal polynomial of $\alpha$, 
$(x-1)^{2n}g\Big(1+\dfrac{2}{x-1}\Big)$ is a minimal polynomial of $\dfrac{y_r}{u_r}$ by Lemma \ref{key2}. By the same lemma, the coefficient of $x^{2n-1}$ in $(x-1)^{2n}g\Big(1+\dfrac{2}{x-1}\Big)$ is
\begin{gather}\label{coe2}
\sum^{2n}_{i=0} \big(-2n+2i\big)\dfrac{a_i}{r^{2i}}.
\end{gather}
Clearly there are infinitely many $r\in\mathbb{Q}$ making \eqref{coe2} nonzero. By Claim \ref{lm1}, we have $\mathbb{Q}\Big(\dfrac{y_r}{u_r}\Big)=\mathbb{Q}\Big(\dfrac{y^2_r}{u^2_r}\Big)$ for those $r$. This completes the proof. 
\end{proof}

\section{Proof of Theorem \ref{thm1}}\label{thms}
The proof of Theorem \ref{thm1} is based on explicit computations. Note that we already know the traces of the first six elements in \eqref{eq.1}.\footnote{That is, $2r=\text{tr}\; \rho=\text{tr}\; \rho',
2s=\text{tr}\; \sigma=\text{tr}\; \sigma'$, and $2t=\text{tr}\; \rho\sigma=\text{tr}\; \rho'\sigma'$.} So, to prove the theorem, it is enough to show the traces of the following elements are contained in $\mathbb{Q}(r,s,t)$:
\begin{equation}\label{eq3}
\rho\rho', \rho\sigma', \sigma\rho', \sigma\sigma', 
\rho\sigma\rho', \rho\sigma\sigma', \rho\rho'\sigma', \sigma\rho'\sigma'.
\end{equation}
Also, by the symmetry between \eqref{ex.1} and \eqref{ex.2}, we can further reduce \eqref{eq3} to the following five elements:
\begin{equation}\label{eq4}
\rho\rho', \rho\sigma', \sigma\sigma', \rho\sigma\rho', \rho\sigma\sigma'.
\end{equation}
Before computing the traces of \eqref{eq4}, we first prove a simple lemma which will be used subsequently in the section. 
When two matrices are represented in terms of the quaternion basis, the trace of their product can be simply computed as the following lemma shows: 
\begin{lemma}\label{comk}
Let
\begin{equation*}
\begin{gathered}
A=a_0\bold{1}+a_1\bold{I}+a_2\bold{(J+K)}+a_3\bold{(J-K)},\\  
B=b_0\bold{1}+b_1\bold{I}+b_2\bold{(J+K)}+b_3\bold{(J-K)}. 
\end{gathered}
\end{equation*}
Then their product $AB$ is of the form 
\allowdisplaybreaks
\begin{equation}\label{no.1}
\begin{gathered}
(a_0b_0+a_1b_1+2a_2b_3+2a_3b_2)\bold{1}+b_0\big(a_1\bold{I}+a_2(\bold{J}+\bold{K})+a_3(\bold{J}-\bold{K})\big)
+a_0\big(b_1\bold{I}+b_2(\bold{J}+\bold{K})+b_3(\bold{J}-\bold{K})\big)\\
+(2a_2b_3-2a_3b_2)\bold{I}+(a_1b_2-a_2b_1)(\bold{J}+\bold{K})+(a_3b_1-a_1b_3)(\bold{J}-\bold{K}),
 \end{gathered}
\end{equation}
and (thus) the trace of $AB$ is 
\begin{gather}\label{no1}
2(a_0b_0+a_1b_1+2a_2b_3+2a_3b_2).
\end{gather}
\end{lemma}
\begin{proof}
Since
\begin{gather*}
A=\begin{pmatrix} 
a_0+a_1 & 2a_2 \\ 
2a_3 & a_0-a_1 
\end{pmatrix}\;\text{and }\;
B=\begin{pmatrix} 
b_0+b_1 & 2b_2 \\ 
2b_3 & b_0-b_1 
\end{pmatrix},
\end{gather*}
$AB$ is equal to 
\allowdisplaybreaks
\begin{align*}
&\begin{pmatrix} 
(a_0+a_1)(b_0+b_1)+4a_2b_3 && (a_0+a_1)2b_2+2a_2(b_0-b_1) \\ 
2a_3(b_0+b_1)+(a_0-a_1)2b_3 && 4a_3b_2+(a_0-a_1)(b_0-b_1) 
\end{pmatrix}\\
\\
=&\begin{pmatrix} 
a_0b_0+a_1b_0+a_0b_1+a_1b_1+4a_2b_3 && 2a_0b_2+2a_1b_2+2a_2b_0-2a_2b_1 \\ 
2a_3b_0+2a_3b_1+2a_0b_3-2a_1b_3 && 4a_3b_2+a_0b_0-a_1b_0-a_0b_1+a_1b_1 
\end{pmatrix}\\
\\
=&\begin{pmatrix} 
a_0b_0+a_1b_0+a_0b_1+a_1b_1+4a_2b_3 && 0 \\ 
0 && 4a_3b_2+a_0b_0-a_1b_0-a_0b_1+a_1b_1 
\end{pmatrix}\\
&+\begin{pmatrix} 
0 && 2a_0b_2+2a_1b_2+2a_2b_0-2a_2b_1 \\ 
2a_3b_0+2a_3b_1+2a_0b_3-2a_1b_3 && 0 
\end{pmatrix}\\
\\
=&\begin{pmatrix} 
a_0b_0+a_1b_1+2a_2b_3 && 0 \\ 
0 && 2a_3b_2+a_0b_0+a_1b_1 
\end{pmatrix}
+\begin{pmatrix} 
a_1b_0+a_0b_1+2a_2b_3 && 0 \\ 
0 && 2a_3b_2-a_1b_0-a_0b_1 
\end{pmatrix}\\
&+\begin{pmatrix} 
0 && 2a_0b_2+2a_1b_2+2a_2b_0-2a_2b_1 \\ 
0 && 0 
\end{pmatrix}
+\begin{pmatrix} 
0 && 0\\ 
2a_3b_0+2a_3b_1+2a_0b_3-2a_1b_3 && 0 
\end{pmatrix}\\
\\
=&\begin{pmatrix} 
a_0b_0+a_1b_1+2a_2b_3+2a_3b_2 && 0 \\ 
0 && 2a_3b_2+a_0b_0+a_1b_1+2a_2b_3 
\end{pmatrix}\\
&+\begin{pmatrix} 
a_1b_0+a_0b_1+2a_2b_3-2a_3b_2 && 0 \\ 
0 && 2a_3b_2-a_1b_0-a_0b_1-2a_2b_3 
\end{pmatrix}\\
&+\begin{pmatrix} 
0 && 2a_0b_2+2a_1b_2+2a_2b_0-2a_2b_1 \\ 
0 && 0 
\end{pmatrix}
+\begin{pmatrix} 
0 && 0\\ 
2a_3b_0+2a_3b_1+2a_0b_3-2a_1b_3 && 0 
\end{pmatrix}\\
\\
=&\begin{pmatrix} 
a_0b_0+a_1b_1+2a_2b_3+2a_3b_2 && 0 \\ 
0 && a_0b_0+a_1b_1+2a_2b_3+2a_3b_2 
\end{pmatrix}\\
&+\begin{pmatrix} 
a_0b_1 && 0 \\ 
0 && -a_0b_1 
\end{pmatrix}
+\begin{pmatrix} 
a_1b_0 && 0 \\ 
0 && -a_1b_0 
\end{pmatrix}
+\begin{pmatrix} 
2a_2b_3-2a_3b_2 && 0 \\ 
0 && 2a_3b_2-2a_2b_3 
\end{pmatrix}\\
&+\begin{pmatrix} 
0 && 2a_0b_2 \\ 
0 && 0 
\end{pmatrix}
+\begin{pmatrix} 
0 && 2a_2b_0 \\ 
0 && 0 
\end{pmatrix}
+\begin{pmatrix} 
0 && 2a_1b_2-2a_2b_1 \\ 
0 && 0 
\end{pmatrix}\\
&+\begin{pmatrix} 
0 && 0\\ 
2a_0b_3 && 0 
\end{pmatrix}
+\begin{pmatrix} 
0 && 0\\ 
2a_3b_0 && 0 
\end{pmatrix}
+\begin{pmatrix} 
0 && 0\\ 
2a_3b_1-2a_1b_3 && 0 
\end{pmatrix}.
\end{align*}
Since
\begin{gather*}
\bold{1}=
\begin{pmatrix} 
1 & 0 \\ 
0 & 1 
\end{pmatrix},
\quad\bold{I}=
\begin{pmatrix} 
1 & 0 \\ 
0 & -1 
\end{pmatrix},
\quad\bold{J}+\bold{K}=
\begin{pmatrix} 
0 & 2 \\ 
0 & 0 
\end{pmatrix},
\text{ and  }\;\;\bold{J}-\bold{K}=
\begin{pmatrix} 
0 & 0 \\ 
2 & 0 
\end{pmatrix},
\end{gather*}
the last given formula above is equal to
\begin{gather*}
(a_0b_0+a_1b_1+2a_2b_3+2a_3b_2)\bold{1}+(a_0b_1)\bold{I}+(a_1b_0)\bold{I}+(2a_2b_3-2a_3b_2)\bold{I}\\
+a_0b_2(\bold{J}+\bold{K})+a_2b_0(\bold{J}+\bold{K})+(a_1b_2-a_2b_1)(\bold{J}+\bold{K})+a_0b_3(\bold{J}-\bold{K})+a_3b_0(\bold{J}-\bold{K})+(a_3b_1-a_1b_3)(\bold{J}-\bold{K}),
\end{gather*}
which is the same as \eqref{no.1}.
\end{proof}

Now we prove Theorem \ref{thm1}. We first compute the traces of the first three elements in \eqref{eq4}. 
\begin{lemma} \label{lem1}
The traces of $\rho\rho'$, $\rho\sigma'$ and $\sigma\sigma'$ are contained in $\mathbb{Q}(r, s, t)$.
\end{lemma}
\begin{proof} 
Recall the formulas of $\rho$ and $\rho'$ in \eqref{ex.1} and \eqref{ex.2}. By Lemma \ref{comk}, the trace of $\rho\rho'$ is equal to
\allowdisplaybreaks
\begin{align*}
&2\Bigg\{r^2+\bigg(\dfrac{r(c+1)}{\sqrt{c^2-1}}\bigg)\bigg(\dfrac{r(c+1)}{\sqrt{c^2-1}}\bigg)+2\bigg(-\dfrac{\tilde{M}\big(2rs-t+(c+\sqrt{c^2-1})t\big)}{2(c-1+2s^2)}\bigg)\bigg(-\dfrac{2rs-t+\big(c-\sqrt{c^2-1}\big)t}{2\tilde{M}(c-1)}\bigg)\\
&+2\bigg(\dfrac{2rs-t+\big(c-\sqrt{c^2-1}\big)t}{2\tilde{M}(c-1)}\bigg)\bigg(\dfrac{\tilde{M}\big(2rs-t+(c+\sqrt{c^2-1})t\big)}{2(c-1+2s^2)}\bigg)\Bigg\}\\
\end{align*}
which can be simplified as
\begin{align*}
&2r^2+\dfrac{2r^2(c+1)^2}{c^2-1}+\dfrac{\big(2rs-t+(c+\sqrt{c^2-1})t\big)\big(2rs-t+\big(c-\sqrt{c^2-1}\big)t\big)}{(c-1+2s^2)(c-1)}\\
&+\dfrac{\big(2rs-t+\big(c-\sqrt{c^2-1}\big)t\big)\big(2rs-t+(c+\sqrt{c^2-1})t\big)}{(c-1)(c-1+2s^2)}\\
=&2r^2+\dfrac{2r^2(c+1)^2}{c^2-1}+\dfrac{\big((2rs-t+tc)+t\sqrt{c^2-1}\big)\big((2rs-t+tc)-t\sqrt{c^2-1}\big)}{(c-1+2s^2)(c-1)}\\
&+\dfrac{\big((2rs-t+tc)-t\sqrt{c^2-1}\big)\big((2rs-t+tc)+t\sqrt{c^2-1}\big)}{(c-1)(c-1+2s^2)}\\
=&2r^2+\dfrac{2r^2(c+1)^2}{c^2-1}+\dfrac{2\big((2rs-t+ct)^2-(c^2-1)t^2\big)}{(c-1+2s^2)(c-1)}.
\end{align*}
Clearly this is contained in $\mathbb{Q}(r,s,t)$. (Recall $c\in \mathbb{Q}(r,s,t)$.)
\\
\\

Similarly, the trace of $\rho\sigma'$ is equal to 
\allowdisplaybreaks
\begin{align*}
2\Bigg\{rs+\bigg(\dfrac{r(c+1)}{\sqrt{c^2-1}}\bigg)\bigg(-\dfrac{s(c+1)}{\sqrt{c^2-1}}\bigg)
+2\bigg(-\dfrac{\tilde{M}\big(2rs-t+(c+\sqrt{c^2-1})t\big)}{2(c-1+2s^2)}\bigg)\bigg(\dfrac{c-1+2s^2}{2\tilde{M}(c-1)}\bigg)\\
+2\bigg(\dfrac{2rs-t+\big(c-\sqrt{c^2-1}\big)t}{2\tilde{M}(c-1)}\bigg)\bigg(-\dfrac{\tilde{M}}{2}\bigg)\Bigg\},
\end{align*}
which can be simplified as  
\allowdisplaybreaks
\begin{align*}
&2rs-\dfrac{2rs(c+1)^2}{c^2-1}
-\dfrac{2rs-t+(c+\sqrt{c^2-1})t}{c-1}
-\dfrac{2rs-t+(c-\sqrt{c^2-1})t}{c-1}\\
=\; &2rs-\dfrac{2rs(c+1)}{c-1}-\dfrac{2(2rs-t+ct)}{c-1}.
\end{align*}
This is also an element of $\mathbb{Q}(r,s,t)$.  
\\
\\

Lastly the trace of $\sigma\sigma'$ is equal to 
\allowdisplaybreaks
\begin{align*}
&2\Bigg\{s^2+\bigg(-\dfrac{s(c+1)}{\sqrt{c^2-1}}\bigg)\bigg(-\dfrac{s(c+1)}{\sqrt{c^2-1}}\bigg)
+2\bigg(\dfrac{\tilde{M}}{2}\bigg)\bigg(\dfrac{c-1+2s^2}{2\tilde{M}(c-1)}\bigg)+2\bigg(-\dfrac{c-1+2s^2}{2\tilde{M}(c-1)}\bigg)\bigg(-\dfrac{\tilde{M}}{2}\bigg)\Bigg\}\\
=\;&2s^2+\dfrac{2s^2(c+1)^2}{c^2-1}+\dfrac{c-1+2s^2}{(c-1)}+\dfrac{c-1+2s^2}{(c-1)}\\
=\;&2s^2+\dfrac{2s^2(c+1)}{c-1}+\dfrac{2(c-1+2s^2)}{(c-1)},
\end{align*}
which is contained in $\mathbb{Q}(r,s,t)$ as well. This completes the proof. 
\end{proof}

Now we compute the traces of $\rho\sigma\rho'$ and $\rho\sigma\sigma'$. 
To do so, we need the explicit formula of $\rho\sigma$ first.
\begin{lemma} \label{lem2}
$\rho\sigma$ is
\begin{align*}
&t\bold{1}+\dfrac{t\sqrt{c^2-1}}{c-1}\bold{I}+\dfrac{(c+1+\sqrt{c^2-1})\tilde{M}\big(r(c-1)+st-st(c+\sqrt{c^2-1})\big)}{2\sqrt{c^2-1}(c-1+2s^2)}(\bold{J}+\bold{K})\\
&+\dfrac{(c+1-\sqrt{c^2-1})\big(r(c-1)+st-st(c-\sqrt{c^2-1})\big)}{2\tilde{M}(c-1)\sqrt{c^2-1}}\big(\bold{J}-\bold{K}\big).
\end{align*}
\end{lemma}
\begin{proof}
We restate the formulas of $\rho$ and $\sigma$ given in \eqref{ex.1}:
\allowdisplaybreaks
\begin{equation*}
\begin{aligned}
&\rho=r\bold{1}+\dfrac{r(c+1)}{\sqrt{c^2-1}}\bold{I}-\dfrac{\tilde{M}\big(2rs-t+(c+\sqrt{c^2-1})t\big)}{2(c-1+2s^2)}(\bold{J}+\bold{K})+\dfrac{2rs-t+\big(c-\sqrt{c^2-1}\big)t}{2\tilde{M}(c-1)}\big(\bold{J}-\bold{K}\big),\\
&\sigma=s\bold{1}-\dfrac{s(c+1)}{\sqrt{c^2-1}}\bold{I}+\dfrac{\tilde{M}}{2}(\bold{J}+\bold{K})-\dfrac{c-1+2s^2}{2\tilde{M}(c-1)}(\bold{J}-\bold{K}).
\end{aligned}
\end{equation*}
To simply the notation, we abbreviate the above $\rho$ and $\sigma$ as  
\allowdisplaybreaks
\begin{equation*}
\begin{aligned}
&\rho=a_0\bold{1}+a_1\bold{I}+a_2(\bold{J}+\bold{K})+a_3(\bold{J}-\bold{K}),\\
&\sigma=b_0\bold{1}+b_1\bold{I}+b_2(\bold{J}+\bold{K})+b_3(\bold{J}-\bold{K}).
\end{aligned}
\end{equation*} 
Note that $a_0b_0+a_1b_1+2a_2b_3+2a_3b_2=t$ (since we got the matrices in \eqref{ex.1} by initially assuming that the trace of $\rho\sigma$ is $2t$). 
To get $\rho\sigma$,  we now compute each term appearing in \eqref{no.1} explicitly as follows: 
\allowdisplaybreaks
\begin{align}
\textbf{(i)}\;&b_0\big(a_1\bold{I}+a_2(\bold{J}+\bold{K})+a_3(\bold{J}-\bold{K})\big)+a_0\big(b_1\bold{I}+b_2(\bold{J}+\bold{K})+b_3(\bold{J}-\bold{K})\big)\notag\\
=&s\bigg(\dfrac{r(c+1)}{\sqrt{c^2-1}}\bold{I}-\dfrac{\tilde{M}\big(2rs-t+(c+\sqrt{c^2-1})t\big)}{2(c-1+2s^2)}(\bold{J}+\bold{K})+\dfrac{2rs-t+\big(c-\sqrt{c^2-1}\big)t}{2\tilde{M}(c-1)}\big(\bold{J}-\bold{K}\big)\bigg)\notag\\
&+r\bigg(-\dfrac{s(c+1)}{\sqrt{c^2-1}}\bold{I}+\dfrac{\tilde{M}}{2}(\bold{J}+\bold{K})-\dfrac{c-1+2s^2}{2\tilde{M}(c-1)}(\bold{J}-\bold{K})\bigg)\notag\\
=&\bigg(\dfrac{sr(c+1)}{\sqrt{c^2-1}}\bold{I}-\dfrac{s\tilde{M}\big(2rs-t+(c+\sqrt{c^2-1})t\big)}{2(c-1+2s^2)}(\bold{J}+\bold{K})+\dfrac{s\big(2rs-t+(c-\sqrt{c^2-1})t\big)}{2\tilde{M}(c-1)}\big(\bold{J}-\bold{K}\big)\bigg)\notag\\
&+\bigg(-\dfrac{rs(c+1)}{\sqrt{c^2-1}}\bold{I}+\dfrac{r\tilde{M}}{2}(\bold{J}+\bold{K})-\dfrac{r(c-1+2s^2)}{2\tilde{M}(c-1)}(\bold{J}-\bold{K})\bigg)\notag\\
=&\bigg(-\dfrac{s\tilde{M}\big(2rs-t+(c+\sqrt{c^2-1})t\big)}{2(c-1+2s^2)}(\bold{J}+\bold{K})+\dfrac{s\big(2rs-t+(c-\sqrt{c^2-1})t\big)}{2\tilde{M}(c-1)}\big(\bold{J}-\bold{K}\big)\bigg)\notag\\
&+\bigg(\dfrac{r\tilde{M}}{2}(\bold{J}+\bold{K})-\dfrac{r(c-1+2s^2)}{2\tilde{M}(c-1)}(\bold{J}-\bold{K})\bigg)\notag\\
=\;&\bigg(\dfrac{r\tilde{M}}{2}-\dfrac{s\tilde{M}\big(2rs-t+(c+\sqrt{c^2-1})t\big)}{2(c-1+2s^2)}\bigg)(\bold{J}+\bold{K})\notag\\
&-\bigg(\dfrac{r(c-1+2s^2)}{2\tilde{M}(c-1)}-\dfrac{s\big(2rs-t+(c-\sqrt{c^2-1})t\big)}{2\tilde{M}(c-1)}\bigg)\big(\bold{J}-\bold{K}\big)\notag\\
=\;&\bigg(\dfrac{r\tilde{M}(c-1+2s^2)}{2(c-1+2s^2)}-\dfrac{s\tilde{M}\big(2rs-t+(c+\sqrt{c^2-1})t\big)}{2(c-1+2s^2)}\bigg)(\bold{J}+\bold{K})\notag\\
&-\bigg(\dfrac{r(c-1+2s^2)}{2\tilde{M}(c-1)}-\dfrac{s\big(2rs-t+(c-\sqrt{c^2-1})t\big)}{2\tilde{M}(c-1)}\bigg)\big(\bold{J}-\bold{K}\big)\notag\\
=\;&\dfrac{r\tilde{M}(c-1+2s^2)-s\tilde{M}\big(2rs-t+(c+\sqrt{c^2-1})t\big)}{2(c-1+2s^2)}(\bold{J}+\bold{K})\notag\\
&-\dfrac{r(c-1+2s^2)-s\big(2rs-t+(c-\sqrt{c^2-1})t\big)}{2\tilde{M}(c-1)}\big(\bold{J}-\bold{K}\big)\notag\\
=\;&\dfrac{\tilde{M}(rc-r+2s^2r)-\tilde{M}\big(2s^2r-st+st(c+\sqrt{c^2-1})\big)}{2(c-1+2s^2)}(\bold{J}+\bold{K})\notag\\
&-\dfrac{rc-r+2s^2r-\big(2s^2r-st+st(c-\sqrt{c^2-1})\big)}{2\tilde{M}(c-1)}\big(\bold{J}-\bold{K}\big)\notag\\
=\;&\dfrac{\tilde{M}\big(r(c-1)+st-st(c+\sqrt{c^2-1})\big)}{2(c-1+2s^2)}(\bold{J}+\bold{K})-\dfrac{r(c-1)+st-st(c-\sqrt{c^2-1})}{2\tilde{M}(c-1)}\big(\bold{J}-\bold{K}\big);\label{no.2}\\
\notag\\
\notag\\
\textbf{(ii)}\; &2(a_2b_3-a_3b_2)\bold{I} \notag\\
=\;&2\bigg\{\bigg(-\dfrac{\tilde{M}\big(2rs-t+(c+\sqrt{c^2-1})t\big)}{2(c-1+2s^2)}\bigg)\bigg(-\dfrac{c-1+2s^2}{2\tilde{M}(c-1)}\bigg)-\bigg(\dfrac{2rs-t+\big(c-\sqrt{c^2-1}\big)t}{2\tilde{M}(c-1)}\bigg)\bigg(\dfrac{\tilde{M}}{2}\bigg)\bigg\}\bold{I}\notag\\
=\;&\bigg(\dfrac{2rs-t+(c+\sqrt{c^2-1})t}{2(c-1)}-\dfrac{2rs-t+(c-\sqrt{c^2-1})t}{2(c-1)}\bigg)\bold{I}\notag\\
=\;&\dfrac{t\sqrt{c^2-1}}{(c-1)}\bold{I}; \label{no.3}\\
\notag\\
\notag\\
\textbf{(iii)}\; &(a_1b_2-a_2b_1)(\bold{J}+\bold{K})\notag\\
=\;&\bigg\{\bigg(\dfrac{r(c+1)}{\sqrt{c^2-1}}\bigg)\bigg(\dfrac{\tilde{M}}{2}\bigg)-\bigg(-\dfrac{\tilde{M}\big(2rs-t+(c+\sqrt{c^2-1})t\big)}{2(c-1+2s^2)}\bigg)\bigg(-\dfrac{s(c+1)}{\sqrt{c^2-1}}\bigg)\bigg\}(\bold{J}+\bold{K})\notag\\
=\;&\bigg(\dfrac{\tilde{M}r(c+1)}{2\sqrt{c^2-1}}-\dfrac{s(c+1)\tilde{M}(2rs-t+(c+\sqrt{c^2-1})t)}{2\sqrt{c^2-1}(c-1+2s^2)}\bigg)(\bold{J}+\bold{K})\notag\\
=\;&\bigg(\dfrac{\tilde{M}r(c+1)(c-1+2s^2)}{2\sqrt{c^2-1}(c-1+2s^2)}-\dfrac{s(c+1)\tilde{M}(2rs-t+(c+\sqrt{c^2-1})t)}{2\sqrt{c^2-1}(c-1+2s^2)}\bigg)(\bold{J}+\bold{K})\notag\\
=\;&\dfrac{(c+1)\tilde{M}\big(r(c-1+2s^2)-2s^2r+st-st(c+\sqrt{c^2-1})\big)}{2\sqrt{c^2-1}(c-1+2s^2)}(\bold{J}+\bold{K})\notag\\
=\;&\dfrac{(c+1)\tilde{M}\big(r(c-1)+st-st(c+\sqrt{c^2-1})\big)}{2\sqrt{c^2-1}(c-1+2s^2)}(\bold{J}+\bold{K});\\
\notag\\
\notag\\
\textbf{(iv)}\; &(a_3b_1-a_1b_3)(\bold{J}-\bold{K})\notag\\
=\;&\bigg\{\bigg(\dfrac{2rs-t+\big(c-\sqrt{c^2-1}\big)t}{2\tilde{M}(c-1)}\bigg)\bigg(-\dfrac{s(c+1)}{\sqrt{c^2-1}}\bigg)-\bigg(\dfrac{r(c+1)}{\sqrt{c^2-1}}\bigg)\bigg(-\dfrac{c-1+2s^2}{2\tilde{M}(c-1)}\bigg)\bigg\}(\bold{J}-\bold{K})\notag\\
=\;&\bigg(-\dfrac{s(c+1)\big(2rs-t+(c-\sqrt{c^2-1})t\big)}{2\tilde{M}(c-1)\sqrt{c^2-1}}+\dfrac{r(c+1)(c-1+2s^2)}{2\tilde{M}(c-1)\sqrt{c^2-1}}\bigg)(\bold{J}-\bold{K})\notag\\
=\;&\dfrac{-s(c+1)\big(2rs-t+(c-\sqrt{c^2-1})t\big)+r(c+1)(c-1+2s^2)}{2\tilde{M}(c-1)\sqrt{c^2-1}}(\bold{J}-\bold{K})\notag\\
=\;&\dfrac{(c+1)\big(-2s^2r+st-st(c-\sqrt{c^2-1})\big)+(c+1)\big(r(c-1)+2s^2r\big)}{2\tilde{M}(c-1)\sqrt{c^2-1}}(\bold{J}-\bold{K})\notag\\
=\;&\dfrac{(c+1)\big(r(c-1)+st-st(c-\sqrt{c^2-1})\big)}{2\tilde{M}(c-1)\sqrt{c^2-1}}(\bold{J}-\bold{K}).\label{no.5}
\end{align}
\\

Combining \eqref{no.2} - \eqref{no.5} together with $a_0b_0+a_1b_1+2a_2b_3+2a_3b_2=t$, we get $\rho\sigma$ is equal to 
\allowdisplaybreaks
\begin{align*}
&t\bold{1}+\dfrac{t\sqrt{c^2-1}}{(c-1)}\bold{I}\\
&+\dfrac{\tilde{M}\big(r(c-1)+st-st(c+\sqrt{c^2-1})\big)}{2(c-1+2s^2)}(\bold{J}+\bold{K})+\dfrac{(c+1)\tilde{M}\big(r(c-1)+st-st(c+\sqrt{c^2-1})\big)}{2\sqrt{c^2-1}(c-1+2s^2)}(\bold{J}+\bold{K})\\
&-\dfrac{\big(r(c-1)+st-st(c-\sqrt{c^2-1})\big)}{2\tilde{M}(c-1)}\big(\bold{J}-\bold{K}\big)+\dfrac{(c+1)\big(r(c-1)+st-st(c-\sqrt{c^2-1})\big)}{2\tilde{M}(c-1)\sqrt{c^2-1}}(\bold{J}-\bold{K})\\
=\;&t\bold{1}+\dfrac{t\sqrt{c^2-1}}{(c-1)}\bold{I}\\
&+\bigg(\dfrac{\sqrt{c^2-1}\tilde{M}\big(r(c-1)+st-st(c+\sqrt{c^2-1})\big)}{2\sqrt{c^2-1}(c-1+2s^2)}+\dfrac{(c+1)\tilde{M}\big(r(c-1)+st-st(c+\sqrt{c^2-1})\big)}{2\sqrt{c^2-1}(c-1+2s^2)}\bigg)(\bold{J}+\bold{K})\\
&+\bigg(-\dfrac{\sqrt{c^2-1}\big(r(c-1)+st-st(c-\sqrt{c^2-1})\big)}{2\tilde{M}(c-1)\sqrt{c^2-1}}+\dfrac{(c+1)\big(r(c-1)+st-st(c-\sqrt{c^2-1})\big)}{2\tilde{M}(c-1)\sqrt{c^2-1}}\bigg)(\bold{J}-\bold{K})\\
=\;&t\bold{1}+\dfrac{t\sqrt{c^2-1}}{c-1}\bold{I}+\dfrac{(c+1+\sqrt{c^2-1})\tilde{M}\big(r(c-1)+st-st(c+\sqrt{c^2-1})\big)}{2\sqrt{c^2-1}(c-1+2s^2)}(\bold{J}+\bold{K})\\
&+\dfrac{(c+1-\sqrt{c^2-1})\big(r(c-1)+st-st(c-\sqrt{c^2-1})\big)}{2\tilde{M}(c-1)\sqrt{c^2-1}}\big(\bold{J}-\bold{K}\big).
\end{align*}
This completes the proof of Lemma \ref{lem2}.
\end{proof}

Using the above lemma, we compute the traces of $\rho\sigma\rho'$ and $\rho\sigma\sigma'$. 
\begin{lemma}\label{lem3}
The trace of $\rho\sigma\rho'$ is equal to
\allowdisplaybreaks
\begin{gather*}
2tr+\dfrac{2tr(c+1)}{c-1}+\dfrac{2(c+1)(r-2st)t}{c-1+2s^2}+\dfrac{2\big(2stc-r(c-1)\big)\big(2rs-t+tc\big)}{(c-1)(c-1+2s^2)}.
\end{gather*}
\end{lemma}

\begin{proof} 
Recall the formulas of $\rho\sigma$ and $\rho'$ given in Lemma \ref{lem2} and \eqref{ex.2}. By Lemma \ref{comk}, the trace of $\sigma\rho\rho'$ is equal to 
\allowdisplaybreaks
\begin{align*}
&2\Bigg\{tr+\bigg(\dfrac{t\sqrt{c^2-1}}{c-1}\bigg)\bigg(\dfrac{r(c+1)}{\sqrt{c^2-1}}\bigg)\notag\\
&+2\bigg(\dfrac{\big(c+1+\sqrt{c^2-1}\big)\tilde{M}\big(r(c-1)+st-st(c+\sqrt{c^2-1})\big)}{2\sqrt{c^2-1}(c-1+2s^2)}\bigg)\bigg(-\dfrac{2rs-t+(c-\sqrt{c^2-1})t}{2\tilde{M}(c-1)}\bigg)\notag\\
&+2\bigg(\dfrac{(c+1-\sqrt{c^2-1})\big(r(c-1)+st-st(c-\sqrt{c^2-1})\big)}{2\tilde{M}(c-1)\sqrt{c^2-1}}\bigg)\bigg(\dfrac{\tilde{M}\big(2rs-t+(c+\sqrt{c^2-1})t\big)}{2(c-1+2s^2)}\bigg)\Bigg\}\notag\\
=&2\Bigg\{tr+\dfrac{tr(c+1)}{c-1}\notag\\
&-\dfrac{\Big(c+1+\sqrt{c^2-1}\Big)\Big(r(c-1)+st-st(c+\sqrt{c^2-1})\Big)\Big(2rs-t+(c-\sqrt{c^2-1})t\Big)}{2(c-1)\sqrt{c^2-1}(c-1+2s^2)}\notag\\
&+\dfrac{\Big(c+1-\sqrt{c^2-1}\Big)\Big(r(c-1)+st-st(c-\sqrt{c^2-1})\Big)\Big(2rs-t+(c+\sqrt{c^2-1})t\Big)}{2(c-1)\sqrt{c^2-1}(c-1+2s^2)}\Bigg\}\notag\\
=&2tr+\dfrac{2tr(c+1)}{c-1}\notag\\
&-\dfrac{\Big(c+1+\sqrt{c^2-1}\Big)\Big(r(c-1)+st-stc-st\sqrt{c^2-1}\Big)\Big(2rs-t+tc-t\sqrt{c^2-1}\Big)}{\sqrt{c^2-1}(c-1)(c-1+2s^2)}\notag\\
&+\dfrac{\Big(c+1-\sqrt{c^2-1}\Big)\Big(r(c-1)+st-stc+st\sqrt{c^2-1}\Big)\Big(2rs-t+tc+t\sqrt{c^2-1}\Big)}{\sqrt{c^2-1}(c-1)(c-1+2s^2)}\notag\\
=&2tr+\dfrac{2tr(c+1)}{c-1}\notag\\
&-\dfrac{\Big(c+1+\sqrt{c^2-1}\Big)\Big((r-st)(c-1)-st\sqrt{c^2-1}\Big)\Big(2rs-t+tc-t\sqrt{c^2-1}\Big)}{\sqrt{c^2-1}(c-1)(c-1+2s^2)}\notag\\
&+\dfrac{\Big(c+1-\sqrt{c^2-1}\Big)\Big((r-st)(c-1)+st\sqrt{c^2-1}\Big)\Big(2rs-t+tc+t\sqrt{c^2-1}\Big)}{\sqrt{c^2-1}(c-1)(c-1+2s^2)}\notag\\
=&2tr+\dfrac{2tr(c+1)}{c-1}\notag\\
&-\dfrac{\Big((c^2-1)(r-st)+\big(-(c+1)st+(r-st)(c-1)\big)\sqrt{c^2-1}-st(c^2-1)\Big)\Big(2rs-t+tc-t\sqrt{c^2-1}\Big)}{\sqrt{c^2-1}(c-1)(c-1+2s^2)}\\
&+\dfrac{\Big((c^2-1)(r-st)+\big((c+1)st-(r-st)(c-1)\big)\sqrt{c^2-1}-st(c^2-1)\Big)\Big(2rs-t+tc+t\sqrt{c^2-1}\Big)}{\sqrt{c^2-1}(c-1)(c-1+2s^2)}\\
=&2tr+\dfrac{2tr(c+1)}{c-1}\notag\\
&-\dfrac{\Big((c^2-1)(r-st)+\big(-2stc+r(c-1)\big)\sqrt{c^2-1}-st(c^2-1)\Big)\Big(2rs-t+tc-t\sqrt{c^2-1}\Big)}{\sqrt{c^2-1}(c-1)(c-1+2s^2)}\\
&+\dfrac{\Big((c^2-1)(r-st)+\big(2stc-r(c-1)\big)\sqrt{c^2-1}-st(c^2-1)\Big)\Big(2rs-t+tc+t\sqrt{c^2-1}\Big)}{\sqrt{c^2-1}(c-1)(c-1+2s^2)}\\
=&2tr+\dfrac{2tr(c+1)}{c-1}\\ 
&-\dfrac{\Big((c^2-1)(r-2st)-\big(2stc-r(c-1)\big)\sqrt{c^2-1}\Big)\Big((2rs-t+tc)-t\sqrt{c^2-1}\Big)}{\sqrt{c^2-1}(c-1)(c-1+2s^2)} \\
&+\dfrac{\Big((c^2-1)(r-2st)+\big(2stc-r(c-1)\big)\sqrt{c^2-1}\Big)\Big((2rs-t+tc)+t\sqrt{c^2-1}\Big)}{\sqrt{c^2-1}(c-1)(c-1+2s^2)}\\
=&2tr+\dfrac{2tr(c+1)}{c-1}+\dfrac{2\big((c^2-1)(r-2st)\big)\big(t\sqrt{c^2-1}\big)}{\sqrt{c^2-1}(c-1)(c-1+2s^2)}+\dfrac{2\Big(\big(2stc-r(c-1)\big)\sqrt{c^2-1}\Big)\Big(2rs-t+tc\Big)}{\sqrt{c^2-1}(c-1)(c-1+2s^2)}\\
=&2tr+\dfrac{2tr(c+1)}{c-1}+\dfrac{2(c+1)(r-2st)t}{c-1+2s^2}+\dfrac{2\big(2stc-r(c-1)\big)\big(2rs-t+tc\big)}{(c-1)(c-1+2s^2)}.
\end{align*}
This completes the proof. 
\end{proof}

\begin{lemma}\label{lem4}
The trace of $\rho\sigma\sigma'$ is equal to 
\begin{align*}
2ts-\dfrac{2ts(c+1)}{c-1}+\dfrac{2\big(r(c-1)+st-stc\big)}{c-1}-\dfrac{2(c+1)st}{c-1}.
\end{align*}
\end{lemma}

\begin{proof}
Recall the formulas of $\rho\sigma$ and $\sigma'$ given in Lemma \ref{lem2} and \eqref{ex.1}. By Lemma \ref{comk}, the trace of $\rho\sigma\sigma'$ is equal to
\allowdisplaybreaks
\begin{align}\label{mel}
&2\Bigg\{ts+\bigg(\dfrac{t\sqrt{c^2-1}}{c-1}\bigg)\bigg(-\dfrac{s(c+1)}{\sqrt{c^2-1}}\bigg)\notag\\
&+2\bigg(\dfrac{\big(c+1+\sqrt{c^2-1}\big)\tilde{M}\big(r(c-1)+st-st(c+\sqrt{c^2-1})\big)}{2\sqrt{c^2-1}(c-1+2s^2)}\bigg)\bigg(\dfrac{c-1+2s^2}{2\tilde{M}(c-1)}\bigg)\notag\\
&+2\bigg(\dfrac{(c+1-\sqrt{c^2-1})\big(r(c-1)+st-st(c-\sqrt{c^2-1})\big)}{2\tilde{M}(c-1)\sqrt{c^2-1}}\bigg)\bigg(-\dfrac{\tilde{M}}{2}\bigg)\Bigg\}.
\end{align} 
This can be simplified as 
\allowdisplaybreaks
\begin{align*}
&2\bigg\{ts-\dfrac{ts(c+1)}{c-1}+\dfrac{(c+1+\sqrt{c^2-1})\big(r(c-1)+st-st(c+\sqrt{c^2-1})\big)}{2\sqrt{c^2-1}(c-1)}\\
&-\dfrac{(c+1-\sqrt{c^2-1})\big(r(c-1)+st-st(c-\sqrt{c^2-1})\big)}{2\sqrt{c^2-1}(c-1)}\bigg\}\\
=&2ts-\dfrac{2ts(c+1)}{c-1}+\dfrac{(c+1+\sqrt{c^2-1})\big(r(c-1)+st-stc-st\sqrt{c^2-1}\big)}{\sqrt{c^2-1}(c-1)}\\
&-\dfrac{(c+1-\sqrt{c^2-1})\big(r(c-1)+st-stc+st\sqrt{c^2-1}\big)}{\sqrt{c^2-1}(c-1)}\\
=&2ts-\dfrac{2ts(c+1)}{c-1}+\dfrac{(c+1+\sqrt{c^2-1})\big(r(c-1)+st-stc-st\sqrt{c^2-1}\big)}{\sqrt{c^2-1}(c-1)}\\
&-\dfrac{(c+1-\sqrt{c^2-1})\big(r(c-1)+st-stc+st\sqrt{c^2-1}\big)}{\sqrt{c^2-1}(c-1)}\\
=&2ts-\dfrac{2ts(c+1)}{c-1}+\dfrac{2\big(r(c-1)+st-stc\big)\sqrt{c^2-1}}{\sqrt{c^2-1}(c-1)}-\dfrac{2(c+1)st\sqrt{c^2-1}}{\sqrt{c^2-1}(c-1)}\\
=&2ts-\dfrac{2ts(c+1)}{c-1}+\dfrac{2\big(r(c-1)+st-stc\big)}{c-1}-\dfrac{2(c+1)st}{c-1},
\end{align*}
which is equal to \eqref{mel}. This completes the proof. 
\end{proof}

By Lemma \ref{lem1}, Lemma \ref{lem3} and Lemma \ref{lem4}, the traces of the elements in \eqref{eq4} are all contained in $\mathbb{Q}(r,s,t)$,  which implies Theorem \ref{thm1}.

\section{Final Remark}\label{Conjecture}
Although Theorem \ref{main1} answers Question \ref{qe1},  we can ask the following stronger question in the spirit of Theorem \ref{KM}:
\begin{question}\label{qe2}
Let $S_g$ ($g \geq 2$) be a closed surface of genus $g$. Given any real number field $K$ and any quaternion algebra $A$ over $K$ such that $A\otimes_K\mathbb{R}\cong M_2(\mathbb{R})$, 
is there a hyperbolic structure on $S_g$ with integral traces such that its invariant trace field and quaternion algebra are equal to $K$ and $A$?
\end{question}

In this case, based on our work, we expect a negative answer. Instead, we propose the following conjecture:
\begin{conjecture}\label{con}
For each $g$ ($g\geq 2$), there exists only a finite number of real number fields and quaternion algebras which arise as the invariant trace field and quaternion algebra of 
a hyperbolic structure of $S_g$ with integral traces.
\end{conjecture}
For each fixed genus surface, as stated above, only finitely many real number fields and quaternion algebras are expected to arise, via integral traces, as invariant trace field and quaternion algebras of hyperbolic structures of it. But,
according to Theorem \ref{KM}, for any real number field $K$ and quaternion algebra $A$ such that $A\otimes_K\mathbb{R}\cong M_2(\mathbb{R})$ and $A\ncong \bigg(\dfrac{1,1}{K}\bigg)$, it is always possible to construct a hyperbolic surface of sufficiently large genus, having $K$ and $A$ as its invariant trace field and quaternion algebra with integral traces. 
Thus, conjecturally, we do not expect any uniform bound on the number of real number fields and quaternion algebras, which are realizable, via integral traces, as the invariants trace fields and quaternion algebras of hyperbolic structures of a fixed genus surface. 
Instead, it would be interesting to find the growth rate of these numbers with respect to genera of surfaces. Lastly we remark, in the spirit of the above conjecture and the discussion so far, both Theorem \ref{KM} and Theorem \ref{main1} are complementary to each other.
\vspace{10 mm}
Department of Mathematics\\
Columbia University\\
2990 Broadway, New York, NY 10027\\
\\
\emph{Email Address}: bogwang.jeon@gmail.com

\begin{thebibliography}{99}


\bibitem{TG} T. ~Gauglhofer, \emph{Trace coordinates of Teichm\"uller spaces of Riemann surfaces}, PhD thesis, EPFL, 2005. 
\bibitem{hun} P. ~Grillet, \emph{Abstract algebra}, Springer, New York, 2007.
\bibitem{K} J. ~Kahn, V. ~Markovic, \emph{Finding cocompact Fuchsian groups of given trace field and quaternion algebra}, Talk at \emph{Geometric structures on $3$-manifolds}, IAS, Oct 2015.
\bibitem{RM} C. ~Maclachlan, A. ~Reid, \emph{The arithmetic of hyperbolic 3-manifolds}, Springer, New York, 2003.
\bibitem{walter} W. ~Neumann, \emph{Realizing arithmetic invariants of hyperbolic 3-manifolds}, Contem. Math. \textbf{541} (2011), 233-246.
\end{thebibliography}
\end{document}